\documentclass[11pt,fleqn]{article}
\usepackage[utf8]{inputenc}
\usepackage[T1]{fontenc}
\usepackage[brazilian,english]{babel}
\usepackage[a4paper,margin=1in]{geometry}
\usepackage{graphicx}
\usepackage{adjustbox}
\graphicspath{{./}{Figuras_V_D/}{Figuras_finais/}{diagnostico/figuras/}}
\usepackage{subfigure}
\usepackage{multirow}
\usepackage{amsfonts}
\usepackage{amsmath}
\usepackage{bm}
\usepackage{booktabs}
\usepackage{comment}
\usepackage{times}
\usepackage{algorithm}
\usepackage{amsmath}

\usepackage{algpseudocode}
\usepackage{setspace}
\usepackage{tcolorbox} 
\usepackage[normalem]{ulem} 
\usepackage{amssymb,amsthm} 
\usepackage{rotating} 
\usepackage{enumerate} 
\usepackage[dvips]{epsfig} 
\usepackage{multicol}
\usepackage{url} 
\usepackage{color} 
\usepackage{xfrac} 
\usepackage{float}  
\usepackage{hyperref} 
\usepackage{orcidlink}

\usepackage{ragged2e}
\usepackage{microtype}
\usepackage{etoolbox}
\usepackage{lettrine}
\ExplSyntaxOn
\seq_new:N \l_article_math_parts_seq
\tl_new:N \l_article_math_line_tl
\tl_new:N \l_article_math_output_tl
\tl_new:N \l_article_math_part_tl
\box_new:N \l_article_math_measure_box
\dim_new:N \l_article_math_width_dim
\NewDocumentCommand{\ColumnMath}{O{0pt}m}
  {
    \group_begin:
    \dim_set:Nn \l_article_math_width_dim { \linewidth - #1 }
    \seq_set_split:Nnn \l_article_math_parts_seq { \MathBreak } {#2}
    \tl_clear:N \l_article_math_line_tl
    \tl_clear:N \l_article_math_output_tl
    \seq_map_inline:Nn \l_article_math_parts_seq
      {
        \tl_set:Nn \l_article_math_part_tl {##1}
        \hbox_set:Nn \l_article_math_measure_box
          { $\displaystyle \tl_use:N \l_article_math_line_tl ##1$ }
        \dim_compare:nNnT
          { \box_wd:N \l_article_math_measure_box }
          > { \l_article_math_width_dim }
          {
            \tl_if_empty:NF \l_article_math_line_tl
              {
                \tl_put_right:NV \l_article_math_output_tl \l_article_math_line_tl
                \tl_put_right:Nn \l_article_math_output_tl { \\ & }
                \tl_clear:N \l_article_math_line_tl
              }
          }
        \tl_if_empty:NT \l_article_math_line_tl
          {
            \tl_if_head_eq_meaning:nNT {##1} \quad
              { \tl_remove_once:Nn \l_article_math_part_tl { \quad } }
            \tl_if_head_eq_meaning:nNT {##1} \qquad
              { \tl_remove_once:Nn \l_article_math_part_tl { \qquad } }
          }
        \tl_put_right:NV \l_article_math_line_tl \l_article_math_part_tl
      }
    \tl_put_right:NV \l_article_math_output_tl \l_article_math_line_tl
    \begin{aligned}[t] & \tl_use:N \l_article_math_output_tl \end{aligned}
    \group_end:
  }
\ExplSyntaxOff

\allowdisplaybreaks[3]
\newcommand{\compactmath}{%
  \fontsize{9}{11}\selectfont
  \setlength{\abovedisplayskip}{0pt}%
  \setlength{\abovedisplayshortskip}{0pt}%
}
\BeforeBeginEnvironment{equation}{\begingroup\compactmath}
\AfterEndEnvironment{equation}{\endgroup}
\BeforeBeginEnvironment{equation*}{\begingroup\compactmath}
\AfterEndEnvironment{equation*}{\endgroup}
\BeforeBeginEnvironment{align}{\begingroup\compactmath}
\AfterEndEnvironment{align}{\endgroup}
\BeforeBeginEnvironment{align*}{\begingroup\compactmath}
\AfterEndEnvironment{align*}{\endgroup}

\usepackage{authblk}

\definecolor{corrblue}{RGB}{0,70,170}
\newcommand{\corr}[1]{{\color{corrblue}#1}}

\definecolor{revred}{RGB}{190,20,20}

\newtheorem{theorem}{Theorem}
\newtheorem{proposition}[theorem]{Proposition}
\newtheorem{lemma}[theorem]{Lemma}
\newtheorem{corollary}[theorem]{Corollary}
\theoremstyle{definition}
\newtheorem{remark}[theorem]{Remark}
\newif\ifproofsinline
\proofsinlinetrue
\newcommand{\TheoryProof}[1]{%
  \ifproofsinline
    \begin{proof}#1\end{proof}%
  \fi}
\newcommand{\AppendixProof}[2]{%
  \begin{proof}[Proof of #1]#2\end{proof}}
\title{Identification of Spatially Varying Thermal Dissipation in Heterogeneous Media Using Integral Transforms and Nonlinear Least Squares}
\author[1, 2]{André J. P. de Oliveira\,\orcidlink{0000-0003-4566-5586}}

\author[3]{Trenton Lau\,\orcidlink{0009-0002-3517-5129}}

\author[2]{Diego C. Knupp\,\orcidlink{0000-0001-9534-5623}\thanks{Corresponding author: \texttt{diegoknupp@iprj.uerj.br}}}

\author[3]{Gary P. T. Choi\,\orcidlink{0000-0001-5407-9111}}

\author[2]{Luiz A. S. Abreu\,\orcidlink{0000-0002-7634-7014}}

\author[2]{Antônio J. Silva Neto\,\orcidlink{0000-0002-9616-6093}}
\affil[1]{Federal Institute of Northern Minas Gerais, IFNMG, Campus Porteirinha, Avenida José Silveira Lopes, nº 429 Vila Serranópolis, Porteirinha - MG, 39520-000, Brazil.}
\affil[2]{Polytechnic Institute, Rio de Janeiro State University, IPRJ/UERJ, Rua Bonfim 25, Vila Amelia, Nova Friburgo, RJ, 28625-570, Brazil.}

\affil[3]{Department of Mathematics, The Chinese University of Hong Kong, Shatin, N.T., Hong Kong.}

\date{}

\begin{document}
\maketitle

\begingroup
\setlength{\parindent}{0pt}
\setlength{\parskip}{0pt}

\noindent
\begin{minipage}[t]{0.27\linewidth}
    \vspace{0pt}
    {\small\scshape Article Info}\par
    \vspace{0.8ex}\hrule\vspace{0.9ex}

    \raggedright
    {\footnotesize\textit{Keywords:}}\par
    {\footnotesize
    Inverse heat transfer problems\\
    Spatially varying thermal dissipation \\ Integral transform technique \\ Spectral truncation regularization \\ Nonlinear least squares\par}
\end{minipage}%
\hfill
\begin{minipage}[t]{0.69\linewidth}
    \vspace{0pt}
    {\small\scshape Abstract}\par
    \vspace{0.8ex}\hrule\vspace{0.9ex}

    \begingroup
    \selectlanguage{english}
    \footnotesize

This work investigates the identification of spatially varying thermal dissipation coefficients in heterogeneous media by combining integral transforms and nonlinear least squares. The unknown coefficient is represented by an eigenfunction expansion constructed from the known thermal properties, resulting in a coupled modal system whose solution depends nonlinearly on the estimated parameters. In this context, a priori and a posteriori strategies for selecting the truncation order based on the discrepancy principle are investigated. The a priori selection directly uses the measured temperatures to determine the problem dimension before optimization, avoiding the need to solve successive inverse problems for different truncation orders. Parameter estimation employs the Gauss--Newton and Levenberg--Marquardt methods, with sensitivities obtained from the modal system, and considers least-squares functionals in the filtered physical and transformed domains. Estimates of the truncation error of the expansion and of the expected residuals of both strategies support the analysis. The methodology is evaluated using synthetic data, including smooth and discontinuous profiles, different thermal excitations, noise levels, and observation resolutions. Across 270 cases, the a priori selection achieved a mean normalized error of $6.99\%$ with a single estimation per case. With a tolerance factor on the discrepancy level, the two a posteriori variants reached $7.45\%$ and $7.74\%$, at computational times 24 and 5 times higher; without it, the errors were $14.38\%$ and $12.34\%$, at times 182 and 162 times higher. The results demonstrate the advantages of prior truncation selection and modal representation in terms of accuracy and computational efficiency for nonlinear identification.

    \par
    \endgroup
\end{minipage}

\endgroup

\twocolumn
\selectlanguage{english}
\justifying

\section{Introduction}
\label{sec:introduction}

Inverse heat transfer problems are typically ill-posed.
Properties, sources, or boundary conditions are recovered from
indirect temperature measurements, and small perturbations in the data
can produce large variations in the estimates. This behavior requires
regularization to control error propagation without eliminating the
useful information contained in the observations \cite{Alifanov1994,Ozisik2021}.

Spectral representations provide a natural framework for formulating
these problems. By expanding the temperature in eigenfunctions, the spatial
field is described by coefficients that capture different scales of the
thermal response. Function identification is then reduced to determining
a finite number of expansion coefficients. In addition to dimensionality
reduction, the modal structure makes it possible to examine which components
of the solution contain sufficient information to support the reconstruction
and which are dominated by noise
\cite{mikhailov1984unified,cotta2016analytical}.

The same basis can be used to represent the unknown parameters
themselves. Spatially varying properties, sources, or dissipation
coefficients are expanded in eigenfunctions, and estimation is performed
on the series coefficients. This parameterization reduces the
function identification problem to a finite-dimensional problem and
allows coupling integrals to be computed in advance.
Different choices of basis, either simplified auxiliary bases or
bases adapted to the known properties, affect modal coupling, the cost of
evaluating the solution, and the sensitivity of temperatures to the parameters
\cite{NAVEIRACOTTA2010,NAVEIRACOTTA2011,NAVEIRACOTTA20111506,Knupp2021}.

The Classical Integral Transform Technique (CITT) and its generalization
(GITT) provide a systematic framework for solving heat transfer
problems, including formulations involving heterogeneous geometries and
properties
\cite{cotta1993,cotta1997heat,cotta2016analytical}.
Projection onto eigenfunctions converts the partial differential equation
into a system of evolution equations for the transformed coefficients.
When the medium coefficients vary spatially, the associated
eigenvalue problem can be solved through expansion in an
auxiliary basis, yielding eigenfunctions adapted to the medium. Contributions
not incorporated into the eigenvalue problem remain in the
transformed system and may introduce coupling between modes.
This structure is particularly convenient for inverse problems,
since the same representation is used for both the model solution
and the parameterization of the unknown function and the treatment of measurements
\cite{cotta1993,cotta1997heat,cotta2016analytical}.

Naveira-Cotta and collaborators established a line of research
combining integral transform solutions, spectral representation
of properties, and inverse analysis. Eigenfunction expansions were
employed to identify spatially varying thermal conductivity
through Bayesian inference \cite{NAVEIRACOTTA2010}, and
subsequently to estimate thermal conductivity and volumetric
heat capacity simultaneously \cite{NAVEIRACOTTA2011}. A relevant
contribution was the use of measured temperatures in the transformed
domain \cite{NAVEIRACOTTA20111506}. Spatial projection enables
data compression and concentration of the information used in
identification, facilitating the efficient treatment of large datasets,
such as those obtained by infrared thermography. This integration was
experimentally demonstrated in the identification of nanocomposite
properties, combining thermographic measurements, data transformation,
GITT solutions, and Bayesian inference \cite{Knupp2012}.

The identification of spatially varying coefficients multiplying
the temperature field has also been investigated in classical
heat transfer models. Trucu et al. \cite{Trucu2008} reconstructed
the blood perfusion coefficient in a transient bioheat model,
combining finite differences and Tikhonov regularization.
Zhuo et al. \cite{ZhuoLesnic2018} investigated the identification of this
coefficient from final-time temperatures or temporal averages,
using gradients obtained from sensitivity and adjoint problems.
Subsequently, Cao and Lesnic \cite{CaoLesnic2021} addressed
the simultaneous identification of a spatially varying reaction
coefficient and a thermal source, employing conjugate gradients
regularized by the discrepancy principle. These studies
highlight the importance of regularization and iterative
procedures in the recovery of distributed coefficients.

Controlling spectral content constitutes another regularization
mechanism in inverse problems \cite{TuanTrong2011, ZhangFuDeng2013, xiong2015spectral}.
Zhang and Yan \cite{ZhangYan2017} investigated
spectral truncation regularization in linear inverse problems,
selecting the truncation order directly from the discrepancy between
the measured data and their eigenfunction projection, with an application
to the identification of a spatial source in a diffusion model.
Knupp and Abreu \cite{Knupp2016}
employed truncated expansions of measured temperatures, with the order
determined by the discrepancy principle, to regularize the
explicit reconstruction of boundary heat flux.
Similar strategies were used in previous studies on source
identification, with \textit{a priori} selection criteria
applied directly to the measurements and
\textit{a posteriori} criteria based on temperatures
computed from the estimated functions
\cite{deOliveira2026explicit,OliveiraNME2026}.
In another approach, Qiao and Xiong \cite{QiaoXiong2023} employed
eigenfunction expansions and iterative Landweber regularization
for the simultaneous reconstruction of a spatial source and the
initial condition in a two-dimensional heat conduction problem,
investigating \textit{a priori} and \textit{a posteriori} criteria
for choosing the regularization parameter. In spectral truncation
regularization, the truncation order controls the spatial resolution
of the reconstruction and the influence of noise. When applied to
the nonlinear identification of distributed parameters, this order
also determines the dimension of the optimization problem, so that
its subsequent selection may require successive estimations.

This issue is particularly relevant to the identification of
spatially varying thermal dissipation coefficients in
heterogeneous media. These coefficients describe distributed
thermal losses and may represent, in reduced models, the effects
of convective boundary conditions present in multidimensional
formulations. Their identification enables the characterization
of local variations in heat transfer mechanisms from
transient temperature measurements. However, their representation
through spectral expansions introduces coupling between modes,
since transforming the product of the unknown coefficient
and the temperature involves both expansions. Consequently,
the truncations of the solution and the parameter jointly influence
the transformed system. Although the coupling matrix depends
linearly on the unknown coefficients, the transient solution
depends nonlinearly on these parameters. In this context, determining
the spectral order before estimation makes it possible to establish
the number of unknowns in advance and avoid successive optimizations
for different truncation orders.

In this work, the identification of a spatially varying
thermal dissipation coefficient is investigated using integral
transforms and nonlinear least squares. The eigenfunctions are
constructed from the known thermal properties,
while the unknown coefficient is represented by an
expansion associated with the same spectral basis as the temperature.
This parameterization allows modal coupling to be expressed through
linear combinations of precomputed matrices and
the same truncation order to be employed in both representations.
\textit{A priori} and \textit{a posteriori} selection strategies
based on the discrepancy principle are investigated,
evaluating reconstruction accuracy and computational cost.
In the \textit{a priori} selection, the order is determined directly
from the measured temperatures, before nonlinear optimization,
whereas the posterior criteria require successive estimations
to evaluate the temperatures predicted by the model.
The a posteriori criteria are examined with and without
a tolerance factor on the discrepancy level. The numerical
study is accompanied by an analysis of the formulation, which
establishes the decay of the truncation error of the expansion
for discontinuous and smooth coefficients, the expected
residuals of the two selection strategies, and conditions for
the local identifiability of the expansion coefficients.
Least-squares
functionals in the filtered physical and transformed domains are
also investigated, including a modal formulation that preserves
the physical-domain metric. Estimation is
performed using Gauss--Newton \cite{NocedalWright2006} and
Levenberg--Marquardt \cite{Levenberg1944,Marquardt1963},
with sensitivities obtained directly from the modal system.
The methodology is evaluated using synthetic data generated
by the finite volume method, considering smooth
and discontinuous coefficients, different thermal excitations, noise
levels, and observation resolutions.

The remainder of this paper is organized as follows.
Section~\ref{sec:problema_direto} presents the thermal model
and its solution by integral transforms.
Section~\ref{sec:problema_inverso} develops the representation
of the unknown coefficient, the treatment of measurements,
the truncation criteria, and the estimation procedures,
together with the associated error and identifiability
results.
Section~\ref{sec:resultados} presents the numerical results
and discusses accuracy, the influence of acquisition conditions,
and computational cost. Finally,
Section~\ref{sec:conclusions} summarizes the conclusions of this work.

\selectlanguage{english}
\section{Formulation and solution of the direct problem}
\label{sec:problema_direto}

Consider heat conduction in a heterogeneous plate of
length $L$ and thickness $l_e$, subjected to a prescribed
heat flux $g(x,t)$ on one of its surfaces, while the
opposite surface exchanges heat with the surroundings at
temperature $T_\infty$. The ends in the longitudinal direction
are thermally insulated. Adopting a lumped formulation across the thickness, the problem is
described by \cite{cotta1997heat, Knupp2021}

\begin{subequations}
\label{eq:modelo_direto}
\begin{align}
&\ColumnMath[2.8em]{%
w(x)\frac{\partial T}{\partial t}\MathBreak
=
\frac{\partial}{\partial x}
\left[
k(x)\frac{\partial T}{\partial x}
\right]\MathBreak
-d(x)\big[T(x,t)-T_\infty\big]\MathBreak
+\frac{g(x,t)}{l_e},%
}
\label{eq:edp_direto}
\\
&\ColumnMath[2.8em]{%
\left.\frac{\partial T}{\partial x}\right|_{x=0}\MathBreak
=
\left.\frac{\partial T}{\partial x}\right|_{x=L}\MathBreak
=0,%
}
\label{eq:contornos_direto}
\\
&\ColumnMath[2.8em]{%
T(x,0)=T_\infty,%
}
\label{eq:inicial_direto}
\end{align}
\end{subequations}
\noindent for $0<x<L$ and $t>0$, where $T(x,t)$
represents the thickness-averaged temperature, while
$k(x)$ and $w(x)$ are the thermal conductivity
and volumetric heat capacity, respectively. The coefficient
$d(x)$ represents thermal dissipation per unit
volume and unit temperature difference. When this
contribution arises exclusively from convection at the
opposite surface, $d(x)=h(x)/l_e$, where $h(x)$ is the
convective heat transfer coefficient.

To develop the solution, the temperature filtered with
respect to the surroundings is defined as

\begin{equation}
\ColumnMath[2.8em]{%
\theta(x,t)\MathBreak
=T(x,t)-T_\infty.%
}
\label{eq:temperatura_filtrada}
\end{equation}

Substituting
Eq.~\eqref{eq:temperatura_filtrada}
into Eqs.~\eqref{eq:modelo_direto} yields

\begin{subequations}
\label{eq:modelo_filtrado}
\begin{align}
&\ColumnMath[2.8em]{%
w(x)\frac{\partial\theta}{\partial t}\MathBreak
=
\frac{\partial}{\partial x}
\left[
k(x)\frac{\partial\theta}{\partial x}
\right]\MathBreak
-d(x)\theta(x,t)\MathBreak
+q_v(x,t),%
}
\label{eq:edp_filtrada}
\\
&\ColumnMath[2.8em]{%
\left.\frac{\partial\theta}{\partial x}\right|_{x=0}\MathBreak
=
\left.\frac{\partial\theta}{\partial x}\right|_{x=L}\MathBreak
=0,%
}
\label{eq:contornos_filtrados}
\\
&\ColumnMath[2.8em]{%
\theta(x,0)=0,%
}
\label{eq:inicial_filtrada}
\end{align}
\end{subequations}

\noindent where

\begin{equation}
\ColumnMath[2.8em]{%
q_v(x,t)\MathBreak
=\frac{g(x,t)}{l_e}%
}
\label{eq:fonte_volumetrica}
\end{equation}
\noindent represents the volumetric source equivalent
to the heat flux applied at the surface.

\subsection{Integral transform and modal system}

To solve Eqs.~\eqref{eq:modelo_filtrado},
the following eigenvalue problem is considered:

\begin{subequations}
\label{eq:autoproblema}
\begin{align}
&\ColumnMath[2.8em]{%
\frac{d}{dx}
\left[
k(x)\frac{d\psi_i(x)}{dx}
\right]\MathBreak
+\mu_i^2w(x)\psi_i(x)\MathBreak
=0,\MathBreak
\qquad 0<x<L,%
}
\\
&\ColumnMath[2.8em]{%
\psi_i'(0)=\psi_i'(L)=0.%
}
\end{align}
\end{subequations}

The coefficient $d(x)$ is excluded from the
eigenvalue problem so that the eigenfunctions depend only
on the known properties $k(x)$ and $w(x)$.

The normalized eigenfunctions are defined by

\begin{equation}
\ColumnMath[2.8em]{%
\widetilde\psi_i(x)\MathBreak
=\frac{\psi_i(x)}{\sqrt{\mathcal N_i}},\MathBreak
\qquad
\mathcal N_i\MathBreak
=\int_0^L w(x)\psi_i^2(x)\,dx,%
}
\label{eq:normalizacao}
\end{equation}
\noindent satisfying the orthogonality property

\begin{equation}
\ColumnMath[2.8em]{%
\int_0^L
w(x)\widetilde\psi_i(x)\widetilde\psi_j(x)\,dx
\MathBreak
=\delta_{ij},%
}
\label{eq:ortogonalidade}
\end{equation}

\noindent where $\delta_{ij}$ is the Kronecker delta.
Ordering the eigenvalues in ascending order, the first
mode corresponds to

\begin{equation}
\ColumnMath[2.8em]{%
\mu_1=0,\MathBreak
\qquad
\widetilde\psi_1(x)=\frac{1}{\sqrt W},\MathBreak
\qquad
W=\int_0^L w(x)\,dx.%
}
\label{eq:modo_constante}
\end{equation}

% ---------------------------------------------------------------------
% [ADDED] Assumptions, weak identity, eigenvalue/eigenfunction bounds
% ---------------------------------------------------------------------
\newcommand{\ProofLemEig}{%
(a) By the min--max principle, $\mu_i^2$ is the $i$-th
min--max value, over $u\in H^1(0,L)$, of the Rayleigh quotient

\begin{equation*}
\ColumnMath[1.5em]{%
R[u]\MathBreak
=\frac{\int_0^Lk(x)\,u'(x)^2\,dx}{\int_0^Lw(x)\,u(x)^2\,dx}.%
}
\end{equation*}

This quotient lies between $(k_{\min}/w_{\max})R_0[u]$ and
$(k_{\max}/w_{\min})R_0[u]$, where $R_0[u]$ is the quotient
obtained with $k\equiv w\equiv1$, whose $i$-th min--max
value is $[(i-1)\pi/L]^2$.

(b) Fix $i\ge2$ and define

\begin{equation*}
\ColumnMath[1.5em]{%
E(x)\MathBreak
=\mu_i^2\,w(x)k(x)\,\widetilde\psi_i(x)^2\MathBreak
+\left[k(x)\widetilde\psi_i'(x)\right]^2.%
}
\end{equation*}

Since $(k\widetilde\psi_i')'=-\mu_i^2w\widetilde\psi_i$,
one obtains $E'=\mu_i^2(wk)'\widetilde\psi_i^2$ almost
everywhere, so that $|E'|\le[|(wk)'|/(wk)]\,E$ and, by
Gronwall's inequality, $E(x)\le e^{\mathcal V}E(y)$ for all
$x,y\in[0,L]$. Integrating in $y$ and using
Eqs.~\eqref{eq:ortogonalidade}
and~\eqref{eq:energia_modal},

\begin{equation*}
\ColumnMath[1.5em]{%
L\,E(x)\MathBreak
\le e^{\mathcal V}\int_0^LE(y)\,dy\MathBreak
\le2k_{\max}e^{\mathcal V}\mu_i^2.%
}
\end{equation*}

Hence $(k\widetilde\psi_i')^2\le E\le C_\psi'^2\mu_i^2$ and
$\widetilde\psi_i^2\le E/(\mu_i^2wk)\le C_\psi^2$.

(c) With two continuous derivatives of $k$ and $w$, the
Liouville transformation reduces Eqs.~\eqref{eq:autoproblema}
to normal form on $(0,\Phi)$ with a bounded potential, for
which the stated behavior holds under second-kind boundary
conditions \cite{Habermann2021,Teschl2012}.%
}

The analysis developed in this work adopts two assumptions on
the coefficients of Eqs.~\eqref{eq:modelo_filtrado}:
(A1) $k$ and $w$ are Lipschitz continuous on $[0,L]$, with
$0<k_{\min}\le k(x)\le k_{\max}$ and
$0<w_{\min}\le w(x)\le w_{\max}$;
(A2) $d$ is bounded, piecewise continuous, non-negative, and $\int_0^L d(x)\,dx>0$ (equivalently, not zero a.e.\ on $(0,L)$).

Multiplying Eqs.~\eqref{eq:autoproblema}, written for
$\widetilde\psi_i$, by $\widetilde\psi_j$, integrating by
parts, and using the boundary conditions together with
Eq.~\eqref{eq:ortogonalidade} gives

\begin{equation}
\ColumnMath[2.8em]{%
\int_0^L
k(x)\widetilde\psi_i'(x)\widetilde\psi_j'(x)\,dx\MathBreak
=\mu_i^2\delta_{ij}.%
}
\label{eq:energia_modal}
\end{equation}

The growth of the eigenvalues and the magnitude of the
eigenfunctions are controlled by the quantities

\begin{equation}
\ColumnMath[2.8em]{%
\Phi\MathBreak
=\int_0^L\sqrt{\frac{w(x)}{k(x)}}\,dx,\MathBreak
\qquad
\mathcal V\MathBreak
=\int_0^L\frac{\left|(wk)'(x)\right|}{w(x)k(x)}\,dx,%
}
\label{eq:phi_v}
\end{equation}

\noindent according to the following result, which is used
in Section~\ref{subsec:representacao_d} to quantify the
truncation error of the expansion of $d(x)$.

\begin{lemma}[Eigenvalue and eigenfunction bounds]
\label{lem:eig}
Under \textup{(A1)}:
\begin{enumerate}
\item[(a)] for every $i\ge1$,

\begin{equation}
\ColumnMath[2.8em]{%
\frac{k_{\min}}{w_{\max}}
\left[\frac{(i-1)\pi}{L}\right]^{2}\MathBreak
\le\mu_i^2\MathBreak
\le\frac{k_{\max}}{w_{\min}}
\left[\frac{(i-1)\pi}{L}\right]^{2};%
}
\label{eq:lim_autovalores}
\end{equation}
\item[(b)] for every $i\ge2$ and every $x\in[0,L]$,

\begin{equation}
\ColumnMath[2.8em]{%
\left|k(x)\widetilde\psi_i'(x)\right|
\le C_\psi'\,\mu_i,\MathBreak
\qquad
\left|\widetilde\psi_i(x)\right|\le C_\psi,%
}
\label{eq:lim_autofuncoes}
\end{equation}
with $C_\psi'=(2k_{\max}e^{\mathcal V}/L)^{1/2}$ and
$C_\psi=C_\psi'/(w_{\min}k_{\min})^{1/2}$, so that the
normalized eigenfunctions are uniformly bounded;
\item[(c)] if, in addition, $k$ and $w$ are twice continuously
differentiable on $[0,L]$, then
$\mu_i=(i-1)\pi/\Phi+O(1/i)$ as $i\to\infty$.
\end{enumerate}
\end{lemma}

\TheoryProof{\ProofLemEig}

Due to the spatial dependence of the properties,
the solution of Eqs.~\eqref{eq:autoproblema} is obtained
through an eigenfunction expansion of an
auxiliary problem with a known solution (see \cite{NAVEIRACOTTA2009}). Using the
cosine eigenfunctions that satisfy the homogeneous
boundary conditions, one writes

\begin{equation}
\ColumnMath[2.8em]{%
\psi_i(x)\MathBreak
\simeq
\sum_{\ell=1}^{M}
b_{\ell i}\varphi_\ell(x),\MathBreak
\qquad
\varphi_\ell(x)
=
\cos\left[
\frac{(\ell-1)\pi x}{L}
\right].%
}
\label{eq:expansao_auxiliar}
\end{equation}

Substituting Eq.~\eqref{eq:expansao_auxiliar}
into Eqs.~\eqref{eq:autoproblema}, projecting onto
the auxiliary basis, and integrating by parts yields
the matrix problem

\begin{equation}
\ColumnMath[2.8em]{%
\mathbf A\mathbf b_i\MathBreak
=\mu_i^2\mathbf B\mathbf b_i,%
}
\label{eq:autoproblema_matricial}
\end{equation}
\noindent where

\begin{subequations}
\label{eq:matrizes_autoproblema}
\begin{align}
&\ColumnMath[2.8em]{%
A_{r\ell}\MathBreak
=
\int_0^L
k(x)\varphi_r'(x)\varphi_\ell'(x)\,dx,%
}
\\
&\ColumnMath[2.8em]{%
B_{r\ell}\MathBreak
=
\int_0^L
w(x)\varphi_r(x)\varphi_\ell(x)\,dx.%
}
\end{align}
\end{subequations}

The auxiliary order $M$ is established from the convergence
of the eigenvalues and eigenfunctions used in the solution
and is distinct from the truncation order of the
temperature expansion.

Using the normalized eigenfunctions, the following
transform--inverse pair is defined \cite{mikhailov1984unified, cotta1993, cotta1997heat, cotta2016analytical}:

\begin{subequations}
\label{eq:transformada}
\begin{align}
&\ColumnMath[2.8em]{%
\bar\theta_i(t)\MathBreak
=
\int_0^L
w(x)\widetilde\psi_i(x)\theta(x,t)\,dx,%
}
\label{eq:transformada_direta}
\\
&\ColumnMath[2.8em]{%
\theta(x,t)\MathBreak
=
\sum_{i=1}^{\infty}
\widetilde\psi_i(x)\bar\theta_i(t).%
}
\label{eq:transformada_inversa}
\end{align}
\end{subequations}

Applying the operator
$\int_0^L\widetilde\psi_i(x)(\cdot)\,dx$
to Eq.~\eqref{eq:edp_filtrada},
integrating the diffusion term by parts, and employing
the boundary conditions and the eigenvalue problem
yields

\begin{equation}
\ColumnMath[2.8em]{%
\frac{d\bar\theta_i}{dt}\MathBreak
+\mu_i^2\bar\theta_i(t)\MathBreak
+\sum_{j=1}^{\infty}
C_{ij}\bar\theta_j(t)\MathBreak
=\bar q_i(t),\MathBreak
\qquad
\bar\theta_i(0)=0,%
}
\label{eq:equacoes_transformadas}
\end{equation}
\noindent where

\begin{subequations}
\label{eq:acoplamento_fonte}
\begin{align}
&\ColumnMath[2.8em]{%
C_{ij}\MathBreak
=
\int_0^L
d(x)\widetilde\psi_i(x)\widetilde\psi_j(x)\,dx,%
}
\\
&\ColumnMath[2.8em]{%
\bar q_i(t)\MathBreak
=
\int_0^L
q_v(x,t)\widetilde\psi_i(x)\,dx\MathBreak
=
\frac{1}{l_e}
\int_0^L
g(x,t)\widetilde\psi_i(x)\,dx.%
}
\end{align}
\end{subequations}

The contribution of $d(x)$ is incorporated into the
coefficients $C_{ij}$, resulting in a coupled system
of ordinary differential equations for the transformed
temperatures.

Truncating the expansion to $N$ terms, the system in
Eq.~\eqref{eq:equacoes_transformadas} is written in
matrix form as

\begin{equation}
\ColumnMath[2.8em]{%
\frac{d\bar{\boldsymbol\theta}}{dt}\MathBreak
+\mathbf H\bar{\boldsymbol\theta}(t)\MathBreak
=\bar{\mathbf q}(t),\MathBreak
\qquad
\bar{\boldsymbol\theta}(0)=\mathbf 0,%
}
\label{eq:sistema_modal}
\end{equation}
\noindent where

\begin{equation}
\ColumnMath[2.8em]{%
\mathbf H\MathBreak
=\boldsymbol\Lambda+\mathbf C,\MathBreak
\qquad
\boldsymbol\Lambda\MathBreak
=
\operatorname{diag}
\left(\mu_1^2,\ldots,\mu_N^2\right),%
}
\label{eq:matriz_sistema}
\end{equation}
\noindent and the vectors
$\bar{\boldsymbol\theta}(t)$ and $\bar{\mathbf q}(t)$
contain the transformed temperatures and sources,
respectively, for the $N$ retained modes.

% ---------------------------------------------------------------------
% [ADDED] Symmetry and positivity of H
% ---------------------------------------------------------------------
\newcommand{\ProofPropH}{%
By Eq.~\eqref{eq:energia_modal},
$\mathbf u^{\mathsf T}\boldsymbol\Lambda\mathbf v
=\int_0^Lk\,u_N'v_N'\,dx$, and by
Eq.~\eqref{eq:acoplamento_fonte},
$\mathbf u^{\mathsf T}\mathbf C\mathbf v
=\int_0^Ld\,u_Nv_N\,dx$, which gives
Eq.~\eqref{eq:forma_H} and the symmetry.
If $\mathbf u^{\mathsf T}\mathbf H\mathbf u=0$, then
$u_N'\equiv0$, so that $u_N$ is a constant $c$ with
$c^2\int_0^Ld\,dx=0$; under (A2), $c=0$ and
$\mathbf u=\mathbf 0$ by the linear independence of the
eigenfunctions. Finally,
$0\le\mathbf u^{\mathsf T}\boldsymbol\Lambda\mathbf u
\le\mu_N^2\|\mathbf u\|_2^2$ and, by
Eq.~\eqref{eq:ortogonalidade},
$\int_0^Ld\,u_N^2\,dx=\int_0^L(d/w)\,w\,u_N^2\,dx$ lies
between $\inf(d/w)\,\|\mathbf u\|_2^2$ and
$\sup(d/w)\,\|\mathbf u\|_2^2$.%
}

The structure of the matrix $\mathbf H$ is established by
the following result.

\begin{proposition}[Symmetry and positivity of $\mathbf H$]
\label{prop:H}
For all $\mathbf u,\mathbf v\in\mathbb R^N$, with
$u_N=\sum_{i=1}^{N}u_i\widetilde\psi_i$ and
$v_N=\sum_{i=1}^{N}v_i\widetilde\psi_i$,

\begin{equation}
\ColumnMath[2.8em]{%
\mathbf u^{\mathsf T}\mathbf H\mathbf v\MathBreak
=\int_0^Lk(x)\,u_N'(x)v_N'(x)\,dx\MathBreak
+\int_0^Ld(x)\,u_N(x)v_N(x)\,dx.%
}
\label{eq:forma_H}
\end{equation}
Consequently, $\mathbf H$ is symmetric and, under
\textup{(A2)}, positive definite, with eigenvalues
satisfying

\begin{equation}
\ColumnMath[2.8em]{%
\inf_{0<x<L}\frac{d(x)}{w(x)}\MathBreak
\le\lambda_{\min}(\mathbf H)\MathBreak
\le\lambda_{\max}(\mathbf H)\MathBreak
\le\mu_N^2+\sup_{0<x<L}\frac{d(x)}{w(x)}.%
}
\label{eq:lim_H}
\end{equation}
\end{proposition}

\TheoryProof{\ProofPropH}

\subsection{Solution of the transformed system}

Since the properties are time-independent,
the matrix $\mathbf H$ is constant. Thus, the solution
of Eq.~\eqref{eq:sistema_modal}, considering the
homogeneous initial condition, is given by

\begin{equation}
\ColumnMath[2.8em]{%
\bar{\boldsymbol\theta}(t)\MathBreak
=
\int_0^t
e^{-\mathbf H(t-\tau)}
\bar{\mathbf q}(\tau)\,d\tau.%
}
\label{eq:solucao_analitica}
\end{equation}

Eq.~\eqref{eq:solucao_analitica} allows the transformed
temperatures to be calculated from the prescribed
source and the system matrix. Since $\mathbf H$
is symmetric and positive definite
(Proposition~\ref{prop:H}), the matrix exponential can be
evaluated through its spectral decomposition, with real
and positive eigenvalues, so that
$\|e^{-\mathbf Ht}\|_2\le1$ for $t\ge0$.

Finally, substituting the transformed coefficients
into Eq.~\eqref{eq:transformada_inversa} and returning
to the definition of the filtered temperature yields

\begin{equation}
\ColumnMath[2.8em]{%
T_N(x,t)\MathBreak
=
T_\infty\MathBreak
+
\sum_{i=1}^{N}
\widetilde\psi_i(x)\bar\theta_i(t).%
}
\label{eq:reconstrucao_temperatura}
\end{equation}

\section{Formulation and solution of the inverse problem}
\label{sec:problema_inverso}

In the inverse problem, the objective is to determine the
thermal dissipation coefficient $d(x)$ from transient
temperature measurements. These measurements are assumed
to be available at positions $x_s$, with $s=1,\ldots,N_x$,
and at times $t_j$, with $j=1,\ldots,N_t$, where
$N_x$ is the number of spatial sensors and $N_t$ is the number
of observation times. The properties $k(x)$ and
$w(x)$, the source $q_v(x,t)$, and the initial and
boundary conditions are assumed to be known.

The filtered measured temperatures are represented by

\begin{equation}
    \ColumnMath[2.8em]{\theta(x_s,t_j)=T(x_s,t_j)-T_\infty},
\end{equation}

\noindent where $\theta(x_s,t_j)$ is the filtered temperature without
measurement errors.

The measurements are then perturbed according to the formulation:

\begin{equation}
\ColumnMath[2.8em]{%
\theta^\delta(x_s,t_j)\MathBreak
=
\theta(x_s,t_j)+\varepsilon_{sj},\MathBreak
\qquad
\varepsilon_{sj}
\sim\operatorname{Normal}(0,\sigma^2),%
}
\label{eq:medidas}
\end{equation}

\noindent where $\varepsilon_{sj}$ represents the error
associated with sensor $s$ at time $t_j$.
The superscript $\delta$ identifies quantities obtained
from data subject to errors. The errors are assumed to
be independent and normally distributed with zero mean
and variance $\sigma^2$, where $\sigma$ is the standard
deviation of the measurement errors.

\subsection{Representation of the dissipation coefficient}
\label{subsec:representacao_d}

To estimate the dissipation coefficient,
the ratio $d(x)/w(x)$ is expanded
in the same eigenfunctions used to represent
the temperature,

\begin{equation}
\ColumnMath[2.8em]{%
\frac{d(x)}{w(x)}\MathBreak
=
\sum_{k=1}^{\infty}
a_k\widetilde\psi_k(x).%
}
\label{eq:expansao_razao}
\end{equation}

Multiplying Eq.~\eqref{eq:expansao_razao} by
$w(x)\widetilde\psi_k(x)$, integrating over the domain,
and employing Eq.~\eqref{eq:ortogonalidade},
the exact expansion coefficients are given by

\begin{equation}
\ColumnMath[2.8em]{%
a_k\MathBreak
=
\int_0^L
d(x)\widetilde\psi_k(x)\,dx.%
}
\label{eq:coeficientes_d}
\end{equation}

Since $d(x)$ is unknown, these coefficients
cannot be calculated directly from
Eq.~\eqref{eq:coeficientes_d} and must be
estimated from the measured temperatures.

Truncating the expansion to $N$ terms and denoting by
$\hat a_k$ the coefficients estimated from the
noisy temperature measurements, the reconstructed
function is written as

\begin{equation}
\ColumnMath[2.8em]{%
\hat d_N(x)\MathBreak
=
w(x)\sum_{k=1}^{N}
\hat a_k\widetilde\psi_k(x).%
}
\label{eq:expansao_d}
\end{equation}

The order $N$ establishes the number of parameters
to be estimated and is determined by the selection
criteria presented below.

Substituting Eq.~\eqref{eq:expansao_d}
into the definition of the coupling coefficients,
Eq.~\eqref{eq:acoplamento_fonte}, yields

\begin{equation}
\ColumnMath[2.8em]{%
\mathbf C(\hat{\mathbf a})\MathBreak
=
\sum_{k=1}^{N}
\hat a_k\mathbf P^{(k)},\MathBreak
\qquad
\hat{\mathbf a}
=
\begin{bmatrix}
\hat a_1&\cdots&\hat a_N
\end{bmatrix}^{\mathsf T},%
}
\label{eq:acoplamento_parametrizado}
\end{equation}

\noindent where $\mathbf P^{(k)}$ are matrices
whose elements correspond to the integrals of
triple products of eigenfunctions,

\begin{equation}
\ColumnMath[2.8em]{%
P_{ij}^{(k)}\MathBreak
=
\int_0^L
w(x)
\widetilde\psi_i(x)
\widetilde\psi_j(x)
\widetilde\psi_k(x)\,dx,\MathBreak
\qquad i,j,k=1,\ldots,N.%
}
\label{eq:produtos_triplos}
\end{equation}

These integrals depend only on the known
eigenfunctions and can be computed before
estimation. Thus, the coupling matrix
corresponding to the reconstructed function is obtained
through linear combinations of the matrices
$\mathbf P^{(k)}$, weighted by the estimated
coefficients $\hat a_k$.

Using the constant mode from
Eq.~\eqref{eq:modo_constante} and the orthogonality
property, it follows that

\begin{equation}
\ColumnMath[2.8em]{%
C_{1j}(\hat{\mathbf a})\MathBreak
=
\frac{\hat a_j}{\sqrt W},\MathBreak
\qquad j=1,\ldots,N.%
}
\label{eq:primeira_linha}
\end{equation}

Each estimated coefficient contributes directly
to the equation for the first mode, associated
with the global energy balance. This structure motivates
the use of the same order $N$ to represent $d(x)$
and the system of transformed temperatures.
The recovery of the coefficients also depends
on the applied excitation and the sensitivity
of the temperatures to the parameters.

% ---------------------------------------------------------------------
% [ADDED] Truncation error of the expansion of d(x)
% ---------------------------------------------------------------------
\newcommand{\ProofPropDecomp}{%
By Eqs.~\eqref{eq:expansao_razao}
and~\eqref{eq:expansao_d},
$(\hat d_N-d)/w
=\sum_{k\le N}(\hat a_k-a_k)\widetilde\psi_k
-\sum_{k>N}a_k\widetilde\psi_k$, and the eigenfunctions
form a complete orthonormal set with respect to the weight
$w(x)$.%
}

\newcommand{\ProofThmStep}{%
(i) From Eqs.~\eqref{eq:autoproblema},
$\widetilde\psi_k=-(k\widetilde\psi_k')'/(\mu_k^2w)$.
For $0\le a<b\le L$, integration by parts gives

\begin{equation*}
\ColumnMath[1.5em]{%
\int_a^b\widetilde\psi_k\,dx\MathBreak
=-\frac{1}{\mu_k^2}
\left[\frac{k\widetilde\psi_k'}{w}\right]_a^b\MathBreak
-\frac{1}{\mu_k^2}
\int_a^bk\widetilde\psi_k'\,\frac{w'}{w^2}\,dx.%
}
\end{equation*}

Since $k\widetilde\psi_k'$ vanishes at $x=0$ and $x=L$,
Lemma~\ref{lem:eig}(b) yields

\begin{equation*}
\ColumnMath[1.5em]{%
\left|\int_{x_\ell}^L\widetilde\psi_k\,dx\right|
\le\frac{C_\psi'(1+\varpi)}{w_{\min}\mu_k},\MathBreak
\qquad
\left|\int_0^L\widetilde\psi_k\,dx\right|
\le\frac{C_\psi'\varpi}{w_{\min}\mu_k}.%
}
\end{equation*}

The bound on $|a_k|$ follows from
$a_k=c_0\int_0^L\widetilde\psi_k\,dx
+\sum_\ell\Delta_\ell\int_{x_\ell}^L\widetilde\psi_k\,dx$.
By Lemma~\ref{lem:eig}(a),
$\mu_k^{-2}\le(w_{\max}/k_{\min})\,L^2/[(k-1)\pi]^2$, and
$\sum_{j\ge n}j^{-2}\le2/n$, so that Eq.~\eqref{eq:cauda}
gives
$\|d-d_n\|_{L^2}^2
\le2w_{\max}^2K_d^2L^2/(k_{\min}\pi^2n)$.

(ii) In the uniform case,
$\widetilde\psi_k=\sqrt{2/(w_0L)}\cos(\nu_kx)$, with
$\nu_k=(k-1)\pi/L$, for $k\ge2$. Since $\sin(\nu_kL)=0$,

\begin{equation*}
\ColumnMath[1.5em]{%
a_k\MathBreak
=-\sqrt{\frac{2}{w_0L}}\,\frac{1}{\nu_k}
\sum_{\ell=1}^{N_d}\Delta_\ell\sin(\nu_kx_\ell),%
}
\end{equation*}

\noindent and, since
$\|w_0\widetilde\psi_k\|_{L^2}^2=w_0$,

\begin{equation*}
\ColumnMath[1.5em]{%
\|d-d_n\|_{L^2}^2\MathBreak
=\frac{2}{L}\sum_{k>n}\frac{1}{\nu_k^2}
\left[\sum_{\ell=1}^{N_d}
\Delta_\ell\sin(\nu_kx_\ell)\right]^2.%
}
\end{equation*}

Expanding the square,

\begin{equation*}
\ColumnMath[1.5em]{%
\left[\sum_{\ell}\Delta_\ell\sin(\nu_kx_\ell)\right]^2\MathBreak
=\frac12\sum_{\ell}\Delta_\ell^2\MathBreak
-\frac12\sum_{\ell}\Delta_\ell^2\cos(2\nu_kx_\ell)\MathBreak
+\sum_{\ell<\ell'}\Delta_\ell\Delta_{\ell'}
\cos\left[\nu_k(x_{\ell'}-x_\ell)\right]\MathBreak
-\sum_{\ell<\ell'}\Delta_\ell\Delta_{\ell'}
\cos\left[\nu_k(x_\ell+x_{\ell'})\right].%
}
\end{equation*}

Every oscillatory term has the form $\cos[(k-1)\pi\xi]$
with $0<\xi<2$, because $0<x_\ell<L$. Its partial sums in
$k$ are therefore bounded by $1/|\sin(\pi\xi/2)|$ and, by
Abel summation against the decreasing weights $(k-1)^{-2}$,
its contribution is $O(n^{-2})$. The constant term gives
$(L/\pi^2)\sum_\ell\Delta_\ell^2\sum_{j\ge n}j^{-2}
=(L/(\pi^2n))\sum_\ell\Delta_\ell^2+O(n^{-2})$.%
}

\newcommand{\ProofPropNeumann}{%
Two integrations by parts, using
$w\widetilde\psi_k=-(k\widetilde\psi_k')'/\mu_k^2$ and
$\widetilde\psi_k'(0)=\widetilde\psi_k'(L)=0$, give

\begin{equation*}
\ColumnMath[1.5em]{%
a_k\MathBreak
=\int_0^L\rho\,w\widetilde\psi_k\,dx\MathBreak
=\frac{1}{\mu_k^2}\int_0^Lk\rho'\,\widetilde\psi_k'\,dx,%
}
\end{equation*}

\noindent and one further integration by parts gives
Eq.~\eqref{eq:ak_suave}. The bound on $|a_k|$ follows from
Lemma~\ref{lem:eig}(b) and from
$\|\widetilde\psi_k\|_{L^2(0,L)}\le w_{\min}^{-1/2}$. By
Lemma~\ref{lem:eig}(a) and
$\sum_{j\ge n}j^{-4}\le\tfrac43n^{-3}$,
Eq.~\eqref{eq:cauda} then gives
Eq.~\eqref{eq:taxa_suave}. In the uniform case,

\begin{equation*}
\ColumnMath[1.5em]{%
\sin\frac{\pi x}{L}\MathBreak
=\frac{2}{\pi}
-\frac{4}{\pi}\sum_{r=1}^{\infty}
\frac{\cos(2r\pi x/L)}{4r^2-1},%
}
\end{equation*}

\noindent and retaining the modes $k\le n$ leaves the terms
with $2r\ge n$, whose mean square is
$\tfrac12\sum_{2r\ge n}16/[\pi^2(4r^2-1)^2]
=4/(3\pi^2n^3)\,[1+o(1)]$.%
}

The representation in Eq.~\eqref{eq:expansao_d} limits the
attainable accuracy through the terms discarded from
Eq.~\eqref{eq:expansao_razao}. To quantify this effect, let
$\rho(x)=d(x)/w(x)$ and let

\begin{equation}
\ColumnMath[2.8em]{%
d_n(x)\MathBreak
=w(x)\sum_{k=1}^{n}a_k\widetilde\psi_k(x)%
}
\label{eq:d_truncado}
\end{equation}

\noindent denote the truncated expansion with the exact
coefficients of Eq.~\eqref{eq:coeficientes_d}.

\begin{proposition}[Error decomposition]
\label{prop:decomp}
For any $\hat{\mathbf a}\in\mathbb R^N$, the function
$\hat d_N$ of Eq.~\eqref{eq:expansao_d} satisfies

\begin{equation}
\ColumnMath[2.8em]{%
\int_0^L\frac{\left[\hat d_N(x)-d(x)\right]^2}{w(x)}\,dx\MathBreak
=\sum_{k=1}^{N}\left(\hat a_k-a_k\right)^2\MathBreak
+\sum_{k>N}a_k^2.%
}
\label{eq:pitagoras}
\end{equation}
\end{proposition}

\TheoryProof{\ProofPropDecomp}

The estimation and truncation errors therefore add in
quadrature, and the last sum in Eq.~\eqref{eq:pitagoras}
cannot be reduced by any estimate of order $N$. Since
$d-d_n=w\,(\rho-\sum_{k\le n}a_k\widetilde\psi_k)$,

\begin{equation}
\ColumnMath[2.8em]{%
\|d-d_n\|_{L^2(0,L)}^2\MathBreak
\le w_{\max}\sum_{k>n}a_k^2.%
}
\label{eq:cauda}
\end{equation}

The decay of this tail depends on the regularity of $d(x)$
and is established below for piecewise-constant and for
smooth coefficients. The results are expressed through the
continuous counterpart of the normalized error adopted in
Section~\ref{sec:resultados},

\begin{equation}
\ColumnMath[2.8em]{%
\operatorname{NRMSE}_{\mathrm c}(n)\MathBreak
=\frac{\|d-d_n\|_{L^2(0,L)}}
{\sqrt L\left[\max d-\min d\right]}.%
}
\label{eq:nrmse_continuo}
\end{equation}

\begin{theorem}[Piecewise-constant coefficient]
\label{thm:step}
Let

\begin{equation}
\ColumnMath[2.8em]{%
d(x)\MathBreak
=c_0+\sum_{\ell=1}^{N_d}
\Delta_\ell\,\mathcal H(x-x_\ell),%
}
\label{eq:d_constante_partes}
\end{equation}
where $\mathcal H$ is the Heaviside step function,
$\Delta_\ell\neq0$, and $0<x_1<\cdots<x_{N_d}<L$.
\begin{enumerate}
\item[(i)] Under \textup{(A1)}, with
$\varpi=(L/w_{\min})\sup|w'|$ and

\begin{equation}
\ColumnMath[2.8em]{%
K_d\MathBreak
=\frac{C_\psi'}{w_{\min}}
\left[(1+\varpi)\sum_{\ell=1}^{N_d}|\Delta_\ell|
+\varpi\,|c_0|\right],%
}
\label{eq:Kd}
\end{equation}
the coefficients satisfy $|a_k|\le K_d/\mu_k$ for every
$k\ge2$, and, for every $n\ge1$,

\begin{equation}
\ColumnMath[2.8em]{%
\|d-d_n\|_{L^2(0,L)}\MathBreak
\le\frac{\sqrt2\,L\,w_{\max}K_d}{\pi\sqrt{k_{\min}}}\;
n^{-1/2}.%
}
\label{eq:taxa_degrau}
\end{equation}
\item[(ii)] In the uniform case, $k\equiv k_0$ and
$w\equiv w_0$,

\begin{equation}
\ColumnMath[2.8em]{%
\|d-d_n\|_{L^2(0,L)}^2\MathBreak
=\frac{L}{\pi^2n}\sum_{\ell=1}^{N_d}\Delta_\ell^2
+O(n^{-2}),%
}
\label{eq:cauda_uniforme}
\end{equation}
so that

\begin{equation}
\ColumnMath[2.8em]{%
\operatorname{NRMSE}_{\mathrm c}(n)\MathBreak
=\frac{C_d}{\sqrt n}\left[1+O(n^{-1})\right],\MathBreak
\qquad
C_d\MathBreak
=\frac{\left(\sum_{\ell}\Delta_\ell^2\right)^{1/2}}
{\pi\left[\max d-\min d\right]}.%
}
\label{eq:Cd}
\end{equation}
In particular, $C_d=1/\pi$ for a single jump and
$C_d=\sqrt2/\pi$ for a rectangular pulse, that is, two
jumps of equal magnitude and opposite signs.
\end{enumerate}
\end{theorem}

\TheoryProof{\ProofThmStep}

\begin{remark}[Heterogeneous medium]
\label{rem:het}
For non-uniform properties, the leading term of the
coefficients is
$a_k\simeq\Delta_\ell\,k(x_\ell)\widetilde\psi_k'(x_\ell)
/[\mu_k^2w(x_\ell)]$ for each jump. Inserting the
leading-order form of the eigenfunctions obtained from the
Liouville transformation gives

\begin{equation}
\ColumnMath[2.8em]{%
C_d^{\mathrm{het}}
=\frac{1}{\pi\left[\max d-\min d\right]}\MathBreak
\times\left[\sum_{\ell=1}^{N_d}
\frac{\Delta_\ell^2\,\Phi/L}
{\sqrt{w(x_\ell)/k(x_\ell)}}\right]^{1/2},%
}
\label{eq:Cd_het}
\end{equation}
so that each jump is weighted by the ratio between the mean
and the local values of $\sqrt{w/k}$. This asymptotic
estimate has the same structure as the rigorous result of
\cite{Habermann2021} for the tail of the eigenfunction
expansion of the Green's function of a regular
Sturm--Liouville problem.
\end{remark}

\begin{proposition}[Smooth coefficient not satisfying the
second-kind boundary conditions]
\label{prop:neumann}
Assume \textup{(A1)} and let $\rho=d/w$ have a
square-integrable second derivative on $(0,L)$. Then, for
every $k\ge2$,

\begin{equation}
\ColumnMath[2.8em]{%
a_k\MathBreak
=\frac{1}{\mu_k^2}
\left[k(x)\,\rho'(x)\,\widetilde\psi_k(x)\right]_{x=0}^{x=L}\MathBreak
-\frac{1}{\mu_k^2}
\int_0^L\left(k\rho'\right)'(x)\,\widetilde\psi_k(x)\,dx.%
}
\label{eq:ak_suave}
\end{equation}
Consequently, $|a_k|\le K_\rho/\mu_k^2$, with

\begin{equation}
\ColumnMath[2.8em]{%
K_\rho\MathBreak
=C_\psi\left[k(0)|\rho'(0)|+k(L)|\rho'(L)|\right]\MathBreak
+\frac{\|(k\rho')'\|_{L^2(0,L)}}{\sqrt{w_{\min}}},%
}
\label{eq:Krho}
\end{equation}
and, for every $n\ge1$,

\begin{equation}
\ColumnMath[2.8em]{%
\|d-d_n\|_{L^2(0,L)}\MathBreak
\le\frac{2\,w_{\max}^{3/2}L^2K_\rho}
{\sqrt3\,\pi^2k_{\min}}\;n^{-3/2}.%
}
\label{eq:taxa_suave}
\end{equation}
A faster decay requires $\rho'(0)=\rho'(L)=0$; otherwise,
the rate $n^{-3/2}$ cannot be improved in general,
regardless of the smoothness of $d$, $k$, and $w$: for
$d(x)=d_0\sin(\pi x/L)$ in a uniform medium,

\begin{equation}
\ColumnMath[2.8em]{%
\operatorname{NRMSE}_{\mathrm c}(n)\MathBreak
=\frac{2}{\sqrt3\,\pi}\,n^{-3/2}\left[1+o(1)\right].%
}
\label{eq:taxa_seno}
\end{equation}
\end{proposition}

\TheoryProof{\ProofPropNeumann}

Theorem~\ref{thm:step} and
Proposition~\ref{prop:neumann} show that the truncation
error decays only algebraically with the order: as
$n^{-1/2}$ for discontinuous coefficients and as $n^{-3/2}$
for smooth coefficients whose ratio $d/w$ has a non-zero
derivative at the ends of the domain, where the
eigenfunctions have zero slope. In the latter case, the
largest deviations of $d_n$ occur near the boundaries.

\subsection{Transformation of measurements}
\label{subsec:transformacao_medidas}

The measured temperatures are projected onto
the same basis used in the solution of the
direct problem. Approximating the transformation
integral by spatial quadrature yields

\begin{equation}
\ColumnMath[2.8em]{%
\bar\theta_i^\delta(t_j)\MathBreak
=
\sum_{s=1}^{N_x}
\omega_s w(x_s)
\widetilde\psi_i(x_s)
\theta^\delta(x_s,t_j),%
}
\label{eq:transformada_medidas}
\end{equation}

\noindent where $\bar\theta_i^\delta(t_j)$
represents the transformed temperature of mode
$i$ at time $t_j$, while $\omega_s$
are the quadrature weights associated with the positions
$x_s$.

For a given order $N$, the retained coefficients
are collected in the vector

\begin{equation}
\ColumnMath[2.8em]{%
\bar{\boldsymbol\theta}_N^\delta(t_j)\MathBreak
=
\begin{bmatrix}
\bar\theta_1^\delta(t_j)&
\cdots&
\bar\theta_N^\delta(t_j)
\end{bmatrix}^{\mathsf T}.%
}
\label{eq:vetor_medidas_transformadas}
\end{equation}

Thus, each measured field, consisting of $N_x$
temperatures, is represented by $N$ transformed
coefficients. The order determines the amount
of modal information used in fitting and
reconstructing the measurements.

\subsection{Solution as a function of the estimated coefficients}
\label{subsec:solucao_parametrizada}

Using the representation of $\hat d_N(x)$
in Eq.~\eqref{eq:expansao_d}, the transformed
system is written as a function of the vector
$\hat{\mathbf a}$ as

\begin{equation}
\ColumnMath[2.8em]{%
\frac{d\bar{\boldsymbol\theta}_N}{dt}\MathBreak
+
\mathbf H_N(\hat{\mathbf a})
\bar{\boldsymbol\theta}_N(t;\hat{\mathbf a})\MathBreak
=
\bar{\mathbf q}_N(t),\MathBreak
\qquad
\bar{\boldsymbol\theta}_N(0;\hat{\mathbf a})
=
\mathbf 0,%
}
\label{eq:sistema_parametrizado}
\end{equation}

\noindent where

\begin{equation}
\ColumnMath[2.8em]{%
\mathbf H_N(\hat{\mathbf a})\MathBreak
=
\boldsymbol\Lambda_N
+
\sum_{k=1}^{N}
\hat a_k\mathbf P^{(k)},\MathBreak
\qquad
\boldsymbol\Lambda_N\MathBreak
=
\operatorname{diag}
\left(\mu_1^2,\ldots,\mu_N^2\right).%
}
\label{eq:matriz_parametrizada}
\end{equation}

The vectors $\bar{\boldsymbol\theta}_N$ and
$\bar{\mathbf q}_N$ contain, respectively,
the transformed temperatures and sources
of the $N$ retained modes.

% ---------------------------------------------------------------------
% [ADDED] Symmetry of H_N(a) at every iterate
% ---------------------------------------------------------------------
\newcommand{\ProofCorHN}{%
By Eqs.~\eqref{eq:expansao_d}
and~\eqref{eq:produtos_triplos},
$\sum_k\hat a_kP_{ij}^{(k)}
=\int_0^L\hat d_N\widetilde\psi_i\widetilde\psi_j\,dx$, so
that Eq.~\eqref{eq:forma_H} holds with $d$ replaced by
$\hat d_N$, and the proof of Proposition~\ref{prop:H}
applies.%
}

\begin{corollary}[Structure of
$\mathbf H_N(\hat{\mathbf a})$]
\label{cor:HN}
For every $\hat{\mathbf a}\in\mathbb R^N$, the matrix
$\mathbf H_N(\hat{\mathbf a})$ is symmetric, and its
eigenvalues lie between $\inf(\hat d_N/w)$ and
$\mu_N^2+\sup(\hat d_N/w)$. It is positive definite
whenever $\hat d_N\ge0$ and $\hat d_N\not\equiv0$.
\end{corollary}

\TheoryProof{\ProofCorHN}

The spectral decomposition is therefore available at every
iterate of the minimization, including iterates for which
$\hat d_N(x)$ changes sign.

Since $\mathbf H_N(\hat{\mathbf a})$
is time-independent, the solution of
Eq.~\eqref{eq:sistema_parametrizado} is given by

\begin{equation}
\ColumnMath[2.8em]{%
\bar{\boldsymbol\theta}_N(t;\hat{\mathbf a})\MathBreak
=
\int_0^t
e^{-\mathbf H_N(\hat{\mathbf a})(t-\tau)}
\bar{\mathbf q}_N(\tau)\,d\tau.%
}
\label{eq:solucao_parametrizada}
\end{equation}

The temperature in the physical domain is then
recovered through the inverse transform,

\begin{equation}
\ColumnMath[2.8em]{%
T_N(x,t;\hat{\mathbf a})\MathBreak
=
T_\infty
+
\sum_{i=1}^{N}
\widetilde\psi_i(x)
\bar\theta_i(t;\hat{\mathbf a}).%
}
\label{eq:temperatura_parametrizada}
\end{equation}

Eqs.~\eqref{eq:solucao_parametrizada}
and~\eqref{eq:temperatura_parametrizada}
allow the temperatures associated with
each estimate of $d(x)$ to be evaluated in
the transformed and physical domains.
Although the coupling matrix is linear
in the coefficients $\hat a_k$, the transient
solution depends on these coefficients through
the matrix exponential. Consequently,
their estimation constitutes a nonlinear
least-squares problem.

The order $N$ remains fixed during each
minimization. It can be determined before
estimation, using the representation of the
measurements, or selected afterward
by evaluating the temperatures calculated
from the estimated functions.

\subsection{Discrepancy principle for truncation order selection}
\label{subsec:criterios_truncamento}

The truncation order determines the number
of modes employed in the representations of
the temperature and the dissipation coefficient.
Low orders may limit the representation
of spatial distributions, whereas high orders
may increase the influence of errors
present in the measurements.
To determine this order, the discrepancy
principle is employed by comparing the
temperature residual with the prescribed
noise variance.

For a filtered temperature field
$\vartheta(x,t)$, the functional is defined as

\begin{equation}
\ColumnMath[2.8em]{%
\mathcal D[\vartheta]\MathBreak
=
\frac{1}{N_xN_t}
\sum_{j=1}^{N_t}
\sum_{s=1}^{N_x}
\left[
\theta^\delta(x_s,t_j)
-
\vartheta(x_s,t_j)
\right]^2,%
}
\label{eq:discrepancia}
\end{equation}

\noindent where $\vartheta$ corresponds
to the representation of the measurements or
the temperature calculated from the estimated
function, depending on the adopted strategy.

Considering a family of fields
$\vartheta_n(x,t)$ associated with the
candidate orders $n$, the smallest order
satisfying the following condition is selected:

\begin{equation}
\ColumnMath[2.8em]{%
N\MathBreak
=
\min\left\{
n\in\{1,\ldots,N_{\max}\}:
\mathcal D[\vartheta_n]\leq\sigma^2
\right\},%
}
\label{eq:criterio_geral_discrepancia}
\end{equation}

\noindent where $\sigma^2$ is the noise
variance and $N_{\max}$ is the maximum
order considered. This limit depends on
the number of available eigenfunctions
and the spatial resolution of the measurements.

% ---------------------------------------------------------------------
% [ADDED] Mean and variance of the discrepancy functional
% ---------------------------------------------------------------------
\newcommand{\ProofLemDisc}{%
By Eqs.~\eqref{eq:medidas} and~\eqref{eq:discrepancia},
the functional is the mean of the $N_xN_t$ independent
terms $(\varepsilon_{sj}+\beta_{sj})^2$. For a Gaussian
variable,
$\mathbb E[(\varepsilon+\beta)^2]=\sigma^2+\beta^2$ and
$\operatorname{Var}[(\varepsilon+\beta)^2]
=2\sigma^4+4\beta^2\sigma^2$, which gives the result.%
}

Since the measurements contain random errors, the functional
in Eq.~\eqref{eq:discrepancia} is a random variable, and its
fluctuation determines the resolution with which the
condition in Eq.~\eqref{eq:criterio_geral_discrepancia} can
be verified.

\begin{lemma}[Mean and variance of the discrepancy]
\label{lem:disc}
Let $\vartheta$ be a deterministic field, independent of the
measurement errors, with
$\beta_{sj}=\theta(x_s,t_j)-\vartheta(x_s,t_j)$ and
$\beta^2=(N_xN_t)^{-1}\sum_{s,j}\beta_{sj}^2$. For the
errors of Eq.~\eqref{eq:medidas},

\begin{equation}
\ColumnMath[2.8em]{%
\mathbb E\left\{\mathcal D[\vartheta]\right\}
=\sigma^2+\beta^2,\MathBreak
\qquad
\operatorname{Var}\left\{\mathcal D[\vartheta]\right\}
=\frac{2\sigma^4+4\sigma^2\beta^2}{N_xN_t}.%
}
\label{eq:media_variancia_discrepancia}
\end{equation}
In particular, for $\beta=0$ the relative standard deviation
of $\mathcal D$ is $\sqrt{2/(N_xN_t)}$.
\end{lemma}

\TheoryProof{\ProofLemDisc}

\subsubsection{A priori selection by truncation of measurements}
\label{subsubsec:selecao_previa}

In the a priori selection, the coefficients obtained
from Eq.~\eqref{eq:transformada_medidas}
are directly employed in the inverse transform.
For each candidate order $n$,
the temperature expanded from the measurements
is defined as

\begin{equation}
\ColumnMath[2.8em]{%
\theta_n^\delta(x,t_j)\MathBreak
=
\sum_{i=1}^{n}
\widetilde\psi_i(x)
\bar\theta_i^\delta(t_j),\MathBreak
\qquad 0\leq x\leq L.%
}
\label{eq:temperatura_expandida_medidas}
\end{equation}

This representation allows the discrepancy
to be evaluated at the sensor positions,

\begin{equation}
\ColumnMath[2.8em]{%
\sigma_n^2\MathBreak
=
\mathcal D[\theta_n^\delta].%
}
\label{eq:discrepancia_previa}
\end{equation}

The order is therefore selected by

\begin{equation}
\ColumnMath[2.8em]{%
N\MathBreak
=
\min\left\{
n\in\{1,\ldots,N_{\max}\}:
\sigma_n^2\leq\sigma^2
\right\}.%
}
\label{eq:selecao_n}
\end{equation}

The procedure uses only the measured
temperatures, the noise level, and the eigenfunctions
calculated from the known properties.
Thus, the order is determined before
estimating $d(x)$, without solving the
transient problem for candidate functions.
After selection, a single estimation
is performed, keeping $N$ fixed in the expansions
of the temperature and the dissipation coefficient.

The use of Eq.~\eqref{eq:selecao_n}
assumes that the discrepancy condition \cite{Mozo1966, deOliveira2026explicit}
is satisfied within the range of orders examined.

\newcommand{\ProofPropApriori}{%
The residual is the mean, over the $N_xN_t$ measurements,
of the squared components of
$(\mathbf I-\boldsymbol\Pi_n)
(\boldsymbol\theta_j+\boldsymbol\varepsilon_j)$,
$j=1,\ldots,N_t$. The cross term has zero mean, and
$\mathbb E\|(\mathbf I-\boldsymbol\Pi_n)
\boldsymbol\varepsilon_j\|_2^2
=\sigma^2\operatorname{tr}[
(\mathbf I-\boldsymbol\Pi_n)^{\mathsf T}
(\mathbf I-\boldsymbol\Pi_n)]$.
Since the rank of $\boldsymbol\Pi_n$ is at most $n$, the
matrix $\mathbf I-\boldsymbol\Pi_n$ acts as the identity on
the null space of $\boldsymbol\Pi_n$, of dimension at least
$N_x-n$. At least $N_x-n$ of its singular values are
therefore not smaller than one, and its squared Frobenius
norm is at least $N_x-n$, with equality exactly for an
orthogonal projector of rank $n$. In the uniform case,
$\boldsymbol\Omega_w=(Lw_0/N_x)\mathbf I$ and, by the
orthogonality of the discrete cosine vectors,
$\boldsymbol\Psi_n^{\mathsf T}\boldsymbol\Psi_n
=[N_x/(Lw_0)]\,\mathbf I_n$ for $n\le N_x$, so that
$\boldsymbol\Pi_n$ is symmetric and idempotent. The residual
is then the orthogonal projection of a Gaussian vector, and
the variance of a non-central chi-square variable with $\nu$
degrees of freedom and non-centrality $\lambda$ is
$2(\nu+2\lambda)$.%
}

To characterize the behavior of $\sigma_n^2$, let
$\boldsymbol\theta_j
=[\theta(x_1,t_j)\ \cdots\ \theta(x_{N_x},t_j)]^{\mathsf T}$
and
$\boldsymbol\theta_j^\delta
=\boldsymbol\theta_j+\boldsymbol\varepsilon_j$ denote the
exact and measured fields sampled at time $t_j$, let
$\boldsymbol\Psi_n$ be the $N_x\times n$ matrix with
elements $\widetilde\psi_i(x_s)$, and define

\begin{equation}
\ColumnMath[2.8em]{%
\boldsymbol\Omega_w\MathBreak
=\operatorname{diag}\left[\omega_sw(x_s)\right],\MathBreak
\qquad
\boldsymbol\Pi_n\MathBreak
=\boldsymbol\Psi_n\boldsymbol\Psi_n^{\mathsf T}
\boldsymbol\Omega_w.%
}
\label{eq:projetor_medidas}
\end{equation}

By Eqs.~\eqref{eq:transformada_medidas}
and~\eqref{eq:temperatura_expandida_medidas}, the field
expanded from the measurements, sampled at the sensor
positions, is $\boldsymbol\Pi_n\boldsymbol\theta_j^\delta$.

\begin{proposition}[Expected a priori residual]
\label{prop:apriori}
The residual of Eq.~\eqref{eq:discrepancia_previa}
satisfies

\begin{equation}
\ColumnMath[2.8em]{%
\mathbb E\left[\sigma_n^2\right]\MathBreak
=\beta_n^2+\kappa_n\sigma^2,\MathBreak
\qquad
\beta_n^2
=\frac{1}{N_xN_t}\sum_{j=1}^{N_t}
\left\|(\mathbf I-\boldsymbol\Pi_n)
\boldsymbol\theta_j\right\|_2^2,%
}
\label{eq:residuo_previo_esperado}
\end{equation}
with

\begin{equation}
\ColumnMath[2.8em]{%
\kappa_n\MathBreak
=\frac{1}{N_x}\operatorname{tr}\left[
(\mathbf I-\boldsymbol\Pi_n)^{\mathsf T}
(\mathbf I-\boldsymbol\Pi_n)\right]\MathBreak
\ge1-\frac{n}{N_x},%
}
\label{eq:kappa}
\end{equation}
and equality holds if and only if $\boldsymbol\Pi_n$ is an
orthogonal projector of rank $n$. This is the case for
uniform $w$, sensors at the midpoints of $N_x$ equal
intervals, weights $\omega_s=L/N_x$, and $n\le N_x$; then
$N_xN_t\,\sigma_n^2/\sigma^2$ has a non-central chi-square
distribution with $(N_x-n)N_t$ degrees of freedom, and
$\operatorname{Var}[\sigma_n^2]
\le2\sigma^4(1+2\beta_n^2/\sigma^2)/(N_xN_t)$.
\end{proposition}

\TheoryProof{\ProofPropApriori}

Each retained mode therefore removes a fraction of
approximately $1/N_x$ of the noise from the residual at
every observation time, and the condition in
Eq.~\eqref{eq:selecao_n} is met, in the mean, at the first
order for which

\begin{equation}
\ColumnMath[2.8em]{%
\beta_n^2\MathBreak
\le(1-\kappa_n)\,\sigma^2\MathBreak
\simeq\frac{n}{N_x}\,\sigma^2.%
}
\label{eq:cruzamento_previo}
\end{equation}

By Lemma~\ref{lem:disc}, this margin exceeds the
fluctuation of the residual whenever
$n>\sqrt{2N_x/N_t}$, so that the selected order is
governed by the decay of $\beta_n^2$, that is, by the
spatial content of the temperature field.

\begin{remark}[Dependence of the selected order on the
noise level]
\label{rem:scaling}
A scaling argument relates the selected order to the noise
level and to the source amplitude. If the
root-mean-square amplitude of the transformed temperatures
decays as $i^{-p}$ and is proportional to a factor
$\alpha$ multiplying the source, then
$\beta_n^2\propto\alpha^2n^{-(2p-1)}$, and
Eq.~\eqref{eq:cruzamento_previo} gives

\begin{equation}
\ColumnMath[2.8em]{%
N\MathBreak
\propto\left(\frac{\alpha\sqrt{N_x}}{\sigma}\right)^{1/p}.%
}
\label{eq:escala_ordem}
\end{equation}
For a discontinuous coefficient,
Theorem~\ref{thm:step} then gives
$\operatorname{NRMSE}\simeq C_d/\sqrt N
\propto(\sigma/\alpha)^{1/(2p)}$ when the truncation error
is dominant. A quasi-static balance of
Eq.~\eqref{eq:equacoes_transformadas} at high order, for a
piecewise-constant source or coefficient, suggests
$p\simeq3$. The order selected a priori thus follows the
spectral content of the temperature field, which is
smoother than $d(x)$, and not that of the coefficient
itself.
\end{remark}

\subsubsection{A posteriori selection}
\label{subsubsec:selecao_posterior}

In the a posteriori selection, the inverse problem
is solved for each candidate order $n$,
employing $n$ terms in the representations
of the temperature and $d(x)$.
The vector estimated at this order is denoted by
$\hat{\mathbf a}^{(n)}$, and the reconstructed
function is

\begin{equation}
\ColumnMath[2.8em]{%
\hat d_n(x)\MathBreak
=
w(x)\sum_{k=1}^{n}
\hat a_k^{(n)}
\widetilde\psi_k(x).%
}
\label{eq:reconstrucao_d_posterior}
\end{equation}

The discrepancy is evaluated using the temperatures
calculated from this function, considering
two criteria.

In the first criterion, denoted P1, the
solution of the fitted model itself is used,
truncated to $n$ modes,

\begin{equation}
\ColumnMath[2.8em]{%
\theta_{\mathrm{P1},n}(x,t)\MathBreak
=
\sum_{i=1}^{n}
\widetilde\psi_i(x)
\bar\theta_i
\left(t;\hat{\mathbf a}^{(n)}\right).%
}
\label{eq:temperatura_p1}
\end{equation}

The transformed coefficients are calculated
from Eq.~\eqref{eq:solucao_parametrizada},
with $N=n$. The corresponding residual is

\begin{equation}
\ColumnMath[2.8em]{%
\sigma_{\mathrm{P1},n}^2\MathBreak
=
\mathcal D[\theta_{\mathrm{P1},n}].%
}
\label{eq:discrepancia_p1}
\end{equation}

In the second criterion, denoted P2, the same
estimated function is inserted into a
direct solution with reference order $N_r$,

\begin{equation}
\ColumnMath[2.8em]{%
\theta_{\mathrm{P2},n}(x,t)\MathBreak
=
\sum_{i=1}^{N_r}
\widetilde\psi_i(x)
\bar\theta_i^{(N_r)}
\left(t;\hat{\mathbf a}^{(n)}\right).%
}
\label{eq:temperatura_p2}
\end{equation}

The superscript $(N_r)$ identifies the order
of the direct system used in this evaluation,
while retaining the representation of $\hat d_n(x)$
with its $n$ coefficients. The corresponding
coupling matrix has elements

\begin{equation}
\ColumnMath[2.8em]{%
C_{ij}^{(N_r)}
\left(\hat{\mathbf a}^{(n)}\right)\MathBreak
=
\int_0^L
\hat d_n(x)
\widetilde\psi_i(x)
\widetilde\psi_j(x)\,dx\MathBreak
=
\sum_{k=1}^{n}
\hat a_k^{(n)}P_{ij}^{(k)},\MathBreak
\qquad i,j=1,\ldots,N_r.%
}
\label{eq:acoplamento_referencia}
\end{equation}

In this expression, the integrals $P_{ij}^{(k)}$
are defined by Eq.~\eqref{eq:produtos_triplos},
with the temperature indices extended
up to $N_r$. The P2 residual is calculated as

\begin{equation}
\ColumnMath[2.8em]{%
\sigma_{\mathrm{P2},n}^2\MathBreak
=
\mathcal D[\theta_{\mathrm{P2},n}].%
}
\label{eq:discrepancia_p2}
\end{equation}

The reference order $N_r$ remains fixed
and is chosen above the candidate orders
employed in the estimation. This procedure
reduces the influence of temperature truncation
on the residual.

The P1 and P2 criteria use the same
estimated function at each candidate order.
The difference lies only in the direct solution
employed to evaluate Eq.~\eqref{eq:discrepancia}.
In both cases, the selected order is

\begin{equation}
\ColumnMath[2.8em]{%
N_{\mathrm P}\MathBreak
=
\min\left\{
n\in\{1,\ldots,N_{\max}\}:
\sigma_{\mathrm P,n}^2\leq\sigma^2
\right\},%
}
\label{eq:selecao_posterior}
\end{equation}

\noindent where $\mathrm P$ identifies
P1 or P2. The noise level and the maximum
candidate order are the same as those used
in the a priori selection, without applying
a tolerance factor to $\sigma^2$.

When no order satisfies the condition in
Eq.~\eqref{eq:selecao_posterior}, the order
whose residual has the smallest
distance from the prescribed variance is used,

\begin{equation}
\ColumnMath[2.8em]{%
N_{\mathrm P}\MathBreak
=
\underset{1\leq n\leq N_{\max}}
{\operatorname{arg\,min}}\,
\left|
\sigma_{\mathrm P,n}^2-\sigma^2
\right|.%
}
\label{eq:selecao_posterior_minimo}
\end{equation}

In the event of a tie, the smallest order
is adopted. This alternative provides an
estimate in cases without a crossing, which
are counted separately in
Section~\ref{subsec:comparacao_truncamento}.

The a posteriori search is performed sequentially,
solving the inverse problem for
$n=1,2,\ldots$ until the condition in
Eq.~\eqref{eq:selecao_posterior} is satisfied.
A total of $N_{\mathrm P}$ estimations
are required when a crossing occurs, and
$N_{\max}$ estimations otherwise. For P2,
each order also requires the evaluation
of the direct solution with $N_r$ modes.
In the a priori selection, the order is
determined from the expansion of the
measurements, and the inverse problem
is solved only once.

\newcommand{\ProofPropAposteriori}{%
The fitted field depends on the data only through
$\hat{\mathbf a}^{(n)}\in\mathbb R^n$, so that
$\mathbf K=(\partial\boldsymbol\vartheta
/\partial\hat{\mathbf a})
(\partial\hat{\mathbf a}/\partial\boldsymbol\varepsilon)$
has rank at most $n$. The linearized residual is
$(\boldsymbol\theta-\boldsymbol\vartheta^\circ)
+(\mathbf I-\mathbf K)\boldsymbol\varepsilon$, whose
expected squared norm is
$N_xN_t\,\beta_{\mathrm P,n}^2
+\sigma^2\|\mathbf I-\mathbf K\|_F^2$, and
$\|\mathbf I-\mathbf K\|_F^2\ge N_xN_t-n$ by the
singular-value argument used in the proof of
Proposition~\ref{prop:apriori}.%
}

\newcommand{\ProofCorCrossing}{%
$(\mathbf I-\mathbf K)\boldsymbol\varepsilon$ is the
orthogonal projection of a Gaussian vector onto a subspace
of dimension $N_xN_t-n$. The normal approximation uses the
mean $\nu$ and the variance $2\nu$ of a chi-square variable
with $\nu$ degrees of freedom, and the last statement
follows from the bound
$\operatorname{Prob}(\chi^2_\nu\ge\nu+2\sqrt{\nu\zeta}
+2\zeta)\le e^{-\zeta}$
\cite{LaurentMassart2000}, with $\nu=N_xN_t-n\le N_xN_t$.%
}

The residuals of the a posteriori criteria behave
differently from $\sigma_n^2$, as established next.

\begin{proposition}[Expected a posteriori residual]
\label{prop:aposteriori}
Fix a candidate order $n$, a criterion
$\mathrm P\in\{\mathrm{P1},\mathrm{P2}\}$, and any of the
least-squares functionals of
Section~\ref{subsec:funcionais}. Let
$\boldsymbol\vartheta(\boldsymbol\varepsilon)
\in\mathbb R^{N_xN_t}$ collect the fitted field of
Eq.~\eqref{eq:temperatura_p1}
or~\eqref{eq:temperatura_p2} at the measurement points, as
a function of the vector $\boldsymbol\varepsilon$ of
measurement errors, and assume that the estimator is
differentiable at $\boldsymbol\varepsilon=\mathbf 0$,
$\boldsymbol\vartheta(\boldsymbol\varepsilon)
=\boldsymbol\vartheta^\circ
+\mathbf K\boldsymbol\varepsilon
+o(\|\boldsymbol\varepsilon\|)$. Then the rank of
$\mathbf K$ is at most $n$ and, within the linearized
residual approximation,

\begin{equation}
\ColumnMath[2.8em]{%
\mathbb E\left[\sigma_{\mathrm P,n}^2\right]\MathBreak
=\beta_{\mathrm P,n}^2
+\frac{\|\mathbf I-\mathbf K\|_F^2}{N_xN_t}\,\sigma^2\MathBreak
\ge\beta_{\mathrm P,n}^2
+\left(1-\frac{n}{N_xN_t}\right)\sigma^2,%
}
\label{eq:residuo_posterior_esperado}
\end{equation}
where
$\beta_{\mathrm P,n}^2
=\|\boldsymbol\theta-\boldsymbol\vartheta^\circ\|_2^2
/(N_xN_t)$ is the misfit of the model fitted to error-free
data and $\boldsymbol\theta$ is the sampled exact field.
\end{proposition}

\TheoryProof{\ProofPropAposteriori}

\begin{corollary}[Probability of satisfying the discrepancy
condition]
\label{cor:crossing}
In the setting of Proposition~\ref{prop:aposteriori},
suppose that the model is exact at order $n$,
$\beta_{\mathrm P,n}=0$, and that $\mathbf K$ is an
orthogonal projector of rank $n$. Then
$N_xN_t\,\sigma_{\mathrm P,n}^2/\sigma^2$ has a chi-square
distribution with $N_xN_t-n$ degrees of freedom and, for
any $\tau\ge1$,

\begin{equation}
\ColumnMath[2.8em]{%
\operatorname{Prob}\left(
\sigma_{\mathrm P,n}^2\le\tau\sigma^2\right)\MathBreak
\approx F_{\mathrm G}\left[
\frac{(\tau-1)N_xN_t+n}{\sqrt{2(N_xN_t-n)}}\right],%
}
\label{eq:prob_cruzamento}
\end{equation}
where $F_{\mathrm G}$ is the standard normal distribution
function. Any misfit $\beta_{\mathrm P,n}>0$ lowers this
probability. Moreover, for every $\zeta>0$, the choice
$\tau\ge1+2\sqrt{\zeta/(N_xN_t)}+2\zeta/(N_xN_t)$ makes it
at least $1-e^{-\zeta}$.
\end{corollary}

\TheoryProof{\ProofCorCrossing}

\begin{remark}[Tolerance factor]
\label{rem:tolerance}
Propositions~\ref{prop:apriori}
and~\ref{prop:aposteriori} show that the two strategies
remove noise from the residual at different rates. In the
a priori selection, the projection acts at every
observation time and the residual decreases by
$\sigma^2/N_x$ per order. In the a posteriori selection,
$n$ parameters are fitted to $N_xN_t$ data and the
decrease is only $\sigma^2/(N_xN_t)$ per order, which is
smaller than the fluctuation
$\sigma^2\sqrt{2/(N_xN_t)}$ of Lemma~\ref{lem:disc}
whenever $n<\sqrt{2N_xN_t}$. With the level $\sigma^2$ of
Eq.~\eqref{eq:selecao_posterior}, the occurrence of a
crossing is then decided essentially by whether the sample
variance of the error realization lies below $\sigma^2$,
and, by Eq.~\eqref{eq:prob_cruzamento} with $\tau=1$, its
probability lies approximately between one half and
$F_{\mathrm G}(1)\approx0.84$ even for an exact model.

For this reason, the discrepancy principle is classically
applied with a factor greater than one
\cite{Mozo1966,BonnefondMarechal2014}, and, for random
errors, residual-based truncation rules are known to fail
or to lose adaptivity when the target order is small
compared with the square root of the number of data
\cite{BlanchardMathe2012,BlanchardHoffmannReiss2018,HarrachJahnPotthast2020,Jahn2021}.
The a posteriori criteria are therefore also considered in
the form

\begin{equation}
\ColumnMath[2.8em]{%
N_{\mathrm P}(\tau)\MathBreak
=\min\left\{
n\in\{1,\ldots,N_{\max}\}:
\sigma_{\mathrm P,n}^2\leq\tau\sigma^2
\right\},%
}
\label{eq:selecao_posterior_tau}
\end{equation}
with $\tau\ge1$, which reduces to
Eq.~\eqref{eq:selecao_posterior} for $\tau=1$. The choice
$\tau=1+3\sqrt{2/(N_xN_t)}$ corresponds to three standard
deviations of the discrepancy and, by
Eq.~\eqref{eq:prob_cruzamento}, makes the crossing of an
exact model practically certain.
\end{remark}

\begin{figure*}[!htbp]
\centering
\includegraphics[width=\textwidth]{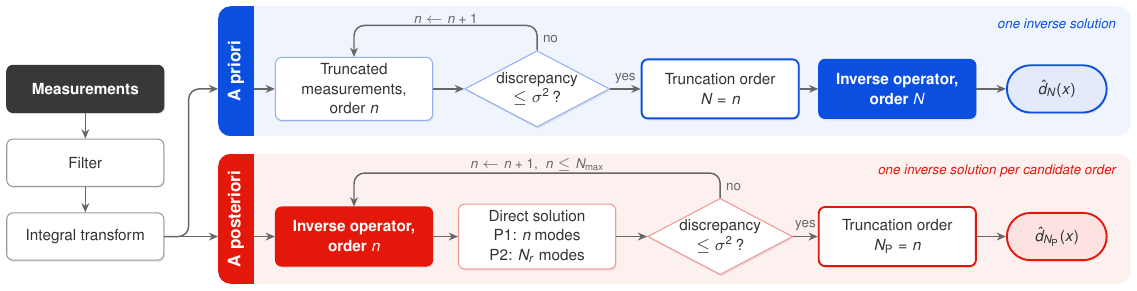}
\caption{Schematic representation of the a priori
truncation order selection procedure based on
the expansion of measurements, and the a posteriori
procedures based on the P1 and P2 criteria.}
\label{fig:flowchart}
\end{figure*}

Figure~\ref{fig:flowchart} presents
the sequence of the a priori and a posteriori
truncation order selection procedures.

\subsection{Least-squares functionals}
\label{subsec:funcionais}

The coefficients are estimated by least
squares, considering fitting in the filtered
physical domain or the transformed domain.
For a given order $N$, the modal residual
vector is defined as

\begin{equation}
\ColumnMath[2.8em]{%
\mathbf e_j(\hat{\mathbf a})\MathBreak
=
\bar{\boldsymbol\theta}_N^\delta(t_j)
-
\bar{\boldsymbol\theta}_N(t_j;\hat{\mathbf a}),%
}
\label{eq:residuo_modal}
\end{equation}

\noindent where the first term contains the
measurement coefficients defined in
Eq.~\eqref{eq:vetor_medidas_transformadas},
while the second is calculated from
Eq.~\eqref{eq:solucao_parametrizada}.
During minimization, $\hat{\mathbf a}$
represents the vector of current estimates.

In formulation A, the modal residuals are
reconstructed at the sensor positions,
resulting in the functional in the filtered
physical domain,

\begin{equation}
\ColumnMath[2.8em]{%
S_A(\hat{\mathbf a})\MathBreak
=
\frac{1}{\sigma^2}
\sum_{j=1}^{N_t}
\mathbf e_j^{\mathsf T}(\hat{\mathbf a})
\boldsymbol\Psi^{\mathsf T}
\boldsymbol\Psi
\mathbf e_j(\hat{\mathbf a}),%
}
\label{eq:funcional_a}
\end{equation}

\noindent where $\boldsymbol\Psi$ is the matrix
of eigenfunction evaluations at the observation
positions, with elements

\begin{equation}
\ColumnMath[2.8em]{%
\Psi_{si}\MathBreak
=
\widetilde\psi_i(x_s),\MathBreak
\qquad s=1,\ldots,N_x,\MathBreak
\qquad i=1,\ldots,N.%
}
\label{eq:matriz_autofuncoes_sensores}
\end{equation}

In formulation B, fitting is performed directly
in the transformed domain,

\begin{equation}
\ColumnMath[2.8em]{%
S_B(\hat{\mathbf a})\MathBreak
=
\frac{1}{\bar\sigma^2}
\sum_{j=1}^{N_t}
\mathbf e_j^{\mathsf T}(\hat{\mathbf a})
\mathbf e_j(\hat{\mathbf a}),%
}
\label{eq:funcional_b}
\end{equation}

\noindent where $\bar\sigma$ is a common scale
for the standard deviations of the transformed errors.

Considering the independent errors in
Eq.~\eqref{eq:medidas}, their propagation through
Eq.~\eqref{eq:transformada_medidas} yields

\begin{equation}
\ColumnMath[2.8em]{%
\bar\sigma_i^2\MathBreak
=
\sigma^2
\sum_{s=1}^{N_x}
\left[
\omega_s w(x_s)\widetilde\psi_i(x_s)
\right]^2,\MathBreak
\qquad
\bar\sigma
=
\frac{1}{N}\sum_{i=1}^{N}\bar\sigma_i.%
}
\label{eq:ruido_modal}
\end{equation}

The transformation may also introduce correlations
between modal errors, whose covariance matrix
has elements

\begin{equation}
\ColumnMath[2.8em]{%
\Sigma_{i\ell}\MathBreak
=
\sigma^2
\sum_{s=1}^{N_x}
\omega_s^2w^2(x_s)
\widetilde\psi_i(x_s)
\widetilde\psi_\ell(x_s).%
}
\label{eq:covariancia_modal}
\end{equation}

In formulation B, only the common scale
from Eq.~\eqref{eq:ruido_modal} is used,
assigning the same weight to all modes.

Formulations A and B employ the same
transformed coefficients, but with different
metrics. Continuous orthogonality with
weight $w(x)$ does not imply that
$\boldsymbol\Psi^{\mathsf T}\boldsymbol\Psi$
is proportional to the identity matrix, so
the minimizers may differ.

To evaluate A directly in modal space,
the following factorization is considered:

\begin{equation}
\ColumnMath[2.8em]{%
\mathbf R^{\mathsf T}\mathbf R\MathBreak
=
\frac{
\boldsymbol\Psi^{\mathsf T}\boldsymbol\Psi
}{\sigma^2}.%
}
\label{eq:fatoracao_metrica}
\end{equation}

Substituting Eq.~\eqref{eq:fatoracao_metrica}
into Eq.~\eqref{eq:funcional_a} yields

\begin{equation}
\ColumnMath[2.8em]{%
S_{A^*}(\hat{\mathbf a})\MathBreak
=
\sum_{j=1}^{N_t}
\left\|
\mathbf R\mathbf e_j(\hat{\mathbf a})
\right\|_2^2\MathBreak
=
S_A(\hat{\mathbf a}).%
}
\label{eq:funcional_a_fatorado}
\end{equation}

For each order considered, $\mathbf R$
is computed before the iterations. This
representation preserves the metric of A
and eliminates the need to reconstruct the residuals
at all sensors during minimization.

The fitting functionals are distinct from
the discrepancy functional in
Eq.~\eqref{eq:discrepancia}, which is used
for truncation selection. In the a posteriori
criteria, the latter is evaluated after
each estimation, using the temperatures
defined for P1 or P2.

\subsection{Sensitivities and estimation procedure}
\label{subsec:sensibilidades_estimacao}

For each order $N$, the coefficients are
determined by minimizing

\begin{equation}
\ColumnMath[2.8em]{%
\underset{\hat{\mathbf a}}{\operatorname{min}}\,
S(\hat{\mathbf a}),%
}
\label{eq:minimizacao_parametros}
\end{equation}

\noindent where $S$ represents $S_A$, $S_B$,
or the equivalent evaluation $S_{A^*}$.
The sensitivities of the transformed temperatures
with respect to the coefficients are defined by

\begin{equation}
\ColumnMath[2.8em]{%
\mathbf s_k(t)\MathBreak
=
\frac{
\partial\bar{\boldsymbol\theta}_N(t;\hat{\mathbf a})
}{
\partial\hat a_k
},\MathBreak
\qquad k=1,\ldots,N.%
}
\label{eq:definicao_sensibilidades}
\end{equation}

Differentiating Eq.~\eqref{eq:sistema_parametrizado}
with respect to $\hat a_k$ yields

\begin{equation}
\ColumnMath[2.8em]{%
\frac{d\mathbf s_k}{dt}\MathBreak
+
\mathbf H_N(\hat{\mathbf a})\mathbf s_k(t)\MathBreak
=
-\mathbf P^{(k)}
\bar{\boldsymbol\theta}_N(t;\hat{\mathbf a}),\MathBreak
\qquad
\mathbf s_k(0)=\mathbf 0,%
}
\label{eq:sensibilidades}
\end{equation}

\noindent since
$\partial\mathbf H_N/\partial\hat a_k
=\mathbf P^{(k)}$.
The system has the same matrix as the
direct problem and admits the representation

\begin{equation}
\ColumnMath[2.8em]{%
\mathbf s_k(t)\MathBreak
=
-\int_0^t
e^{-\mathbf H_N(\hat{\mathbf a})(t-\tau)}
\mathbf P^{(k)}
\bar{\boldsymbol\theta}_N(\tau;\hat{\mathbf a})
\,d\tau.%
}
\label{eq:solucao_sensibilidades}
\end{equation}

To construct the residuals and the Jacobian,
the weighting matrix is defined as

\begin{equation}
\ColumnMath[2.8em]{%
\mathbf L\MathBreak
=
\begin{cases}
\boldsymbol\Psi/\sigma,
& \text{formulation A},\\
\mathbf I_N/\bar\sigma,
& \text{formulation B},\\
\mathbf R,
& \text{formulation A$^*$},
\end{cases}%
}
\label{eq:ponderacao_residuos}
\end{equation}

\noindent where $\mathbf I_N$ is the identity matrix
of order $N$.
The residual blocks and their derivatives are
therefore given by

\begin{equation}
\ColumnMath[2.8em]{%
\mathbf r_j(\hat{\mathbf a})\MathBreak
=
\mathbf L\mathbf e_j(\hat{\mathbf a}),\MathBreak
\qquad
\frac{\partial\mathbf r_j}{\partial\hat a_k}
\MathBreak
=
-\mathbf L\mathbf s_k(t_j).%
}
\label{eq:residuos_sensibilidades}
\end{equation}

The negative sign follows from the definition
of the residual as the difference between the
measurements and the calculated temperatures.

Combining the blocks for all observation times
yields the residual vector $\mathbf r$ and
the Jacobian $\mathbf J$, with

\begin{equation}
\ColumnMath[2.8em]{%
S(\hat{\mathbf a})\MathBreak
=
\mathbf r^{\mathsf T}(\hat{\mathbf a})
\mathbf r(\hat{\mathbf a}),\MathBreak
\qquad
J_{\ell k}\MathBreak
=
\frac{\partial r_\ell}{\partial\hat a_k},%
}
\label{eq:jacobiano}
\end{equation}

\noindent where $\ell$ identifies each
component of the residual vector.

\newcommand{\ProofThmIdent}{%
Let $\mathbf v\in\mathbb R^N$ with
$\mathbf v^{\mathsf T}\mathbf G_c\mathbf v
=\int_0^{t_f}\|\mathbf L\mathbf u(t)\|_2^2\,dt=0$, where
$\mathbf u=\mathbf S\mathbf v$. By (I2) and continuity,
$\mathbf u\equiv\mathbf 0$. By
Eq.~\eqref{eq:sensibilidades} and linearity,
$d\mathbf u/dt+\mathbf H_N(\hat{\mathbf a})\mathbf u
=-\mathbf P(\mathbf v)\bar{\boldsymbol\theta}_N$, with
$\mathbf P(\mathbf v)=\sum_kv_k\mathbf P^{(k)}$, so that
$\mathbf P(\mathbf v)\bar{\boldsymbol\theta}_N(t)
=\mathbf 0$ for all $t$. By
Eq.~\eqref{eq:primeira_linha}, the first row of
$\mathbf P(\mathbf v)$ is
$\mathbf v^{\mathsf T}/\sqrt W$; hence
$\mathbf v^{\mathsf T}\bar{\boldsymbol\theta}_N(t)=0$ for
all $t$, and (I1) gives $\mathbf v=\mathbf 0$. For equally
spaced observation times,
$[t_f/(N_t-1)]\,\mathbf J^{\mathsf T}\mathbf J$ is a
Riemann sum of the continuous integrand and converges to
$\mathbf G_c$ as $N_t\to\infty$; positive definiteness is
preserved for sufficiently large $N_t$ because the smallest
eigenvalue depends continuously on the matrix.%
}

The conditions under which the Jacobian has full rank are
established next. Let
$\mathbf S(t)
=[\,\mathbf s_1(t)\ \cdots\ \mathbf s_N(t)\,]$ collect the
sensitivities of Eq.~\eqref{eq:definicao_sensibilidades},
so that

\begin{equation}
\ColumnMath[2.8em]{%
\mathbf J^{\mathsf T}\mathbf J\MathBreak
=\sum_{j=1}^{N_t}
\mathbf S^{\mathsf T}(t_j)\,
\mathbf L^{\mathsf T}\mathbf L\,
\mathbf S(t_j),%
}
\label{eq:gram_discreto}
\end{equation}

\noindent and define its continuous-time counterpart,

\begin{equation}
\ColumnMath[2.8em]{%
\mathbf G_c(\hat{\mathbf a})\MathBreak
=\int_0^{t_f}
\mathbf S^{\mathsf T}(t)\,
\mathbf L^{\mathsf T}\mathbf L\,
\mathbf S(t)\,dt.%
}
\label{eq:gram_continuo}
\end{equation}

\begin{theorem}[Local identifiability]
\label{thm:ident}
Fix the order $N$, the weighting matrix $\mathbf L$, and
$\hat{\mathbf a}$, and assume that
\begin{enumerate}
\item[(I1)] the vectors
$\bar{\boldsymbol\theta}_N(t;\hat{\mathbf a})$,
$0\le t\le t_f$, span $\mathbb R^N$;
\item[(I2)] $\mathbf L$ has rank $N$.
\end{enumerate}
Then $\mathbf G_c(\hat{\mathbf a})$ is positive definite.
Consequently, for equally spaced observation times there
exists $N_t^\ast$ such that
$\mathbf J^{\mathsf T}\mathbf J$ is positive definite for
all $N_t\ge N_t^\ast$; the Gauss--Newton increment of
Eq.~\eqref{eq:gauss_newton} is then uniquely defined, and a
minimizer with zero residual is isolated.
\end{theorem}

\TheoryProof{\ProofThmIdent}

\begin{remark}[On the hypotheses of
Theorem~\ref{thm:ident}]
\label{rem:ident}
(a) The proof uses only the first row of the coupling
matrix, Eq.~\eqref{eq:primeira_linha}, through which every
coefficient $\hat a_k$ enters the equation of the first
mode.
(b) For separable sources, $q_v(x,t)=Q(x)F(t)$ with
 $F\in L^1(0,t_f)$ and $F\neq0$ a.e.\ on $(0,t_f)$,
condition (I1) holds if and only if
$\mathbf H_N(\hat{\mathbf a})$ has simple eigenvalues and
the vector with components
$\int_0^LQ\widetilde\psi_i\,dx$, $i=1,\ldots,N$, has a
non-zero component along each of its eigenvectors. Indeed,
with the eigenpairs $(\lambda_i,\mathbf z_i)$ of the
symmetric matrix $\mathbf H_N(\hat{\mathbf a})$
(Corollary~\ref{cor:HN}), the transformed temperature is a
combination of the vectors $\mathbf z_i$ with the factors
$\chi_i(t)=\int_0^te^{-\lambda_i(t-s)}F(s)\,ds$, which are
linearly independent on $[0,t_f]$ for distinct
$\lambda_i$, by Titchmarsh's convolution theorem.
(c) In formulation B,
$\mathbf L=\mathbf I_N/\bar\sigma$ and (I2) is automatic;
in formulations A and A$^*$, (I2) requires
$\boldsymbol\Psi$ to have rank $N$, and hence
$N\le N_x$.
(d) The result is local: it characterizes the minimizer in
a neighborhood of $\hat{\mathbf a}$.
\end{remark}

In the Gauss--Newton method, the increment
at iteration $m$ is determined by

\begin{equation}
\ColumnMath[2.8em]{%
\Delta\hat{\mathbf a}_m\MathBreak
=
\underset{\Delta\hat{\mathbf a}}
{\operatorname{arg\,min}}\,
\left\|
\mathbf J_m\Delta\hat{\mathbf a}
+\mathbf r_m
\right\|_2^2,%
}
\label{eq:gauss_newton}
\end{equation}

\noindent where
$\mathbf J_m=\mathbf J(\hat{\mathbf a}_m)$
and $\mathbf r_m=\mathbf r(\hat{\mathbf a}_m)$.
The update is performed by

\begin{equation}
\ColumnMath[2.8em]{%
\hat{\mathbf a}_{m+1}\MathBreak
=
\hat{\mathbf a}_m
+\gamma_m\Delta\hat{\mathbf a}_m,\MathBreak
\qquad 0<\gamma_m\leq1,%
}
\label{eq:atualizacao_gn}
\end{equation}

\noindent where $\gamma_m$ is determined
by a line search using the Armijo criterion.

In the Levenberg--Marquardt method, the following
system is solved:

\begin{equation}
\ColumnMath[2.8em]{%
\left[
\mathbf J_m^{\mathsf T}\mathbf J_m
+
\lambda_{\mathrm{LM},m}
\operatorname{diag}
\left(
\mathbf J_m^{\mathsf T}\mathbf J_m
\right)
\right]
\Delta\hat{\mathbf a}_m\MathBreak
=
-\mathbf J_m^{\mathsf T}\mathbf r_m,%
}
\label{eq:lm}
\end{equation}

\noindent where $\lambda_{\mathrm{LM},m}$
is the damping parameter, adjusted
according to the reduction in the functional.
Accepted increments are used in the update
$\hat{\mathbf a}_{m+1}
=\hat{\mathbf a}_m+\Delta\hat{\mathbf a}_m$.

Minimization is terminated when the functional
stabilizes over consecutive iterations or
when the gradient norm decreases relative
to its initial value, according to the
established tolerances.

In the a priori selection, minimization is performed
with the previously determined order.
In the a posteriori selection, the same procedure
is applied to each candidate order, retaining
the vector corresponding to the order selected
by P1 or P2. The final function is reconstructed
using Eq.~\eqref{eq:expansao_d}, with the
coefficients obtained at the end of minimization.

\begin{remark}[Propagation of measurement errors]
\label{rem:noise}
To first order in $\sigma$, in formulation B and around the
minimizer $\hat{\mathbf a}^\circ$ obtained with error-free
data, the covariance of the estimated coefficients is

\begin{equation}
\ColumnMath[2.8em]{%
\operatorname{Cov}(\hat{\mathbf a})\MathBreak
\simeq\frac{1}{\bar\sigma^{4}}\,
(\mathbf J^{\mathsf T}\mathbf J)^{-1}\,
\mathbf V\,
(\mathbf J^{\mathsf T}\mathbf J)^{-1},\MathBreak
\qquad
\mathbf V
=\sum_{j=1}^{N_t}
\mathbf S^{\mathsf T}(t_j)\,
\boldsymbol\Sigma\,
\mathbf S(t_j),%
}
\label{eq:covariancia_estimativas}
\end{equation}
where $\boldsymbol\Sigma$ is the matrix with elements given
by Eq.~\eqref{eq:covariancia_modal}. When $\boldsymbol\Sigma=\bar\sigma^2\mathbf I_N$ (independent modal data), which holds for uniform $w$ with the midpoint rule and $N\le N_x$, this expression reduces to $(\mathbf J^{\mathsf T}\mathbf J)^{-1}$, so that

\begin{equation}
\ColumnMath[2.8em]{%
\mathbb E\left\|
\hat{\mathbf a}-\hat{\mathbf a}^\circ\right\|_2^2\MathBreak
\simeq\operatorname{tr}\left[
(\mathbf J^{\mathsf T}\mathbf J)^{-1}\right]\MathBreak
\le\frac{N}{\lambda_{\min}
(\mathbf J^{\mathsf T}\mathbf J)},%
}
\label{eq:erro_coeficientes}
\end{equation}
and
$\|\hat{\mathbf a}-\hat{\mathbf a}^\circ\|_2
\le(\sqrt N+\sqrt{2\zeta})
/\lambda_{\min}(\mathbf J^{\mathsf T}\mathbf J)^{1/2}$
with probability at least $1-e^{-\zeta}$
\cite{LaurentMassart2000}. Only the projection of the
errors onto the $N$-dimensional range of $\mathbf J$
contributes to the estimates. Since, by
Eq.~\eqref{eq:gram_discreto},
$\mathbf J^{\mathsf T}\mathbf J$ grows in proportion to
$N_t$ for a fixed observation interval, this contribution
decreases as $N_t^{-1/2}$, and, by
Proposition~\ref{prop:decomp}, it adds in quadrature to
the truncation error. Formulation B employs a common scale
instead of $\boldsymbol\Sigma$ and coincides with the
minimum-variance estimator when
$\boldsymbol\Sigma=\bar\sigma^2\mathbf I_N$. The relation
between a full-rank sensitivity matrix, the Fisher
information matrix, and the standard errors of the
estimates is discussed in \cite{CintronArias2009}, and the
convergence of the Gauss--Newton and Levenberg--Marquardt
iterations for zero- and small-residual problems in
\cite{MadsenNielsenTingleff2004,BergouDiouaneKungurtsev2020}.
\end{remark}

\section{Results and discussion}
\label{sec:resultados}

This section presents the numerical results obtained
from simulated data, evaluating reconstruction
accuracy, the influence of acquisition conditions,
and the computational cost of the methodology.
The simulations were performed on a computer
equipped with an Intel Core i5-10500H processor
operating at $2.50\,\mathrm{GHz}$, $8\,\mathrm{GB}$
of RAM, and a $64$-bit operating system
with an x64 architecture.
The implementations were developed in
Wolfram Mathematica 15.0~\cite{Mathematica15}
and Python 3.14~\cite{Python314}, employing
the Mathematica functions \texttt{NDSolve},
\texttt{DSolve}, and \texttt{NMinimize} for
the numerical and symbolic solution of equations
and for comparing the minimization procedures.

\subsection{Inverse problem parameters and volumetric sources}
\label{subsec:parametros_numericos}

The results are presented for a heterogeneous
plate whose properties vary smoothly along
the longitudinal direction. The thermal conductivity
and volumetric heat capacity are described by

\begin{equation}
\ColumnMath[2.8em]{%
k(x)\MathBreak
=k_0+(k_L-k_0)\mathcal T(x),\MathBreak
\qquad
w(x)
=w_0+(w_L-w_0)\mathcal T(x),%
}
\label{eq:propriedades_teste}
\end{equation}

\noindent where

\begin{equation}
\ColumnMath[2.8em]{%
\mathcal T (x)\MathBreak
=
\left\{
1+\exp\left[-200(x-0.04)\right]
\right\}^{-1},%
}
\label{eq:transicao_propriedades}
\end{equation}

\noindent with $x$ in meters.

Table~\ref{tab:properties} summarizes the geometric and thermophysical parameters used to generate the synthetic data.

\begin{table*}[!tp]
\centering
\small
\caption{Geometry and thermophysical property
parameters of the plate, defined by
Eqs.~\eqref{eq:propriedades_teste}
and~\eqref{eq:transicao_propriedades}.}
\label{tab:properties}
\setlength{\tabcolsep}{7pt}
\renewcommand{\arraystretch}{1.15}
\begin{adjustbox}{max width=\linewidth}
\begin{tabular}{llcc}
\toprule
\textbf{Quantity}
& \textbf{Symbol}
& \textbf{Value}
& \textbf{Unit}
\\
\midrule
Length
& $L$
& $0.08$
& $\mathrm{m}$
\\
Thickness
& $l_e$
& $0.002$
& $\mathrm{m}$
\\
Final time
& $t_f$
& $100$
& $\mathrm{s}$
\\
Lower reference thermal conductivity
& $k_0$
& $0.549$
& $\mathrm{W\,m^{-1}\,K^{-1}}$
\\
Upper reference thermal conductivity
& $k_L$
& $2.29$
& $\mathrm{W\,m^{-1}\,K^{-1}}$
\\
Lower reference volumetric heat capacity
& $w_0$
& $2.23\times10^6$
& $\mathrm{J\,m^{-3}\,K^{-1}}$
\\
Upper reference volumetric heat capacity
& $w_L$
& $2.582\times10^6$
& $\mathrm{J\,m^{-3}\,K^{-1}}$
\\
Ambient temperature
& $T_\infty$
& $25$
& $^\circ\mathrm C$
\\
\bottomrule
\end{tabular}
\end{adjustbox}
\end{table*}

Three profiles are considered for $d(x)$,
allowing the reconstruction of a smooth
function and functions with discontinuities
to be evaluated.
The first has a sinusoidal distribution,

\begin{equation}
\ColumnMath[2.8em]{%
d(x)\MathBreak
=
5\times10^4
\sin\left(\frac{\pi x}{L}\right).%
}
\label{eq:d_seno}
\end{equation}

The second exhibits a step transition,

\begin{equation}
\ColumnMath[2.8em]{%
d(x)\MathBreak
=
1.5\times10^4
+3\times10^4\mathcal H(x-0.035),%
}
\label{eq:d_degrau}
\end{equation}

\noindent while the third corresponds
to a spatial pulse,

\begin{equation}
\ColumnMath[2.8em]{%
d(x)\MathBreak
=
10^4
+4\times10^4
\mathcal I_{[0.025,\,0.055]}(x).%
}
\label{eq:d_pulso}
\end{equation}

In Eqs.~\eqref{eq:d_seno}--\eqref{eq:d_pulso},
$d(x)$ is expressed in
$\mathrm{W\,m^{-3}\,K^{-1}}$ and $x$ in meters.
The function $\mathcal H$ represents the
Heaviside step function, with $\mathcal H(0)=1$,
and the indicator function is defined by

\begin{equation}
\ColumnMath[2.8em]{%
\mathcal I_{[a,b]}(x)\MathBreak
=
\begin{cases}
1, & a\leq x\leq b,\\
0, & \text{otherwise}.
\end{cases}%
}
\label{eq:funcao_indicadora}
\end{equation}

To evaluate the influence of heating
on estimation \cite{Knupp2021}, three
distributions are considered for the volumetric
source $q_v(x,t)$. The first is spatially uniform
and sinusoidal in time,

\begin{equation}
\ColumnMath[2.8em]{%
q_1(x,t)\MathBreak
=
5\times10^5
\sin\left(\frac{\pi t}{50}\right).%
}
\label{eq:fonte_q1}
\end{equation}

The second acts in four zones, with
constant intensity over time,

\begin{equation}
\ColumnMath[2.8em]{%
q_2(x,t)\MathBreak
=
10^6\mathcal Z(x),%
}
\label{eq:fonte_q2}
\end{equation}

\noindent and the third combines this
spatial distribution with the temporal
dependence of $q_1$,

\begin{equation}
\ColumnMath[2.8em]{%
q_3(x,t)\MathBreak
=
10^6\mathcal Z(x)
\sin\left(\frac{\pi t}{50}\right).%
}
\label{eq:fonte_q3}
\end{equation}

The sources are expressed in
$\mathrm{W\,m^{-3}}$, with $t$ in seconds.
The active zones are defined by

\begin{equation}
\ColumnMath[2.8em]{%
\mathcal Z(x)\MathBreak
=
\sum_{j=1}^{4}
\mathcal I_{
\left[
(4j-3)L/16,\,
(4j-1)L/16
\right]
}(x).%
}
\label{eq:zonas_ativas}
\end{equation}

These zones have a width of
$1\,\mathrm{cm}$, are centered at $x=1$,
$3$, $5$, and $7\,\mathrm{cm}$, and occupy half
of the plate length.
The volumetric amplitudes are obtained
by dividing the surface heat fluxes
by the thickness $l_e$.
Since the active zones occupy half
of the domain, $q_1$ and $q_3$ have the same
source integrated along $x$ at each
instant. The source $q_2$ maintains this
integrated intensity at its maximum value
throughout the experiment. For $q_1$ and $q_3$, the sinusoidal dependence
completes one period over the interval considered,
with heat supplied between $0$ and
$50\,\mathrm{s}$ and extracted between $50$
and $100\,\mathrm{s}$.

\begin{figure*}[!htbp]
\centering
\includegraphics[width=\textwidth]{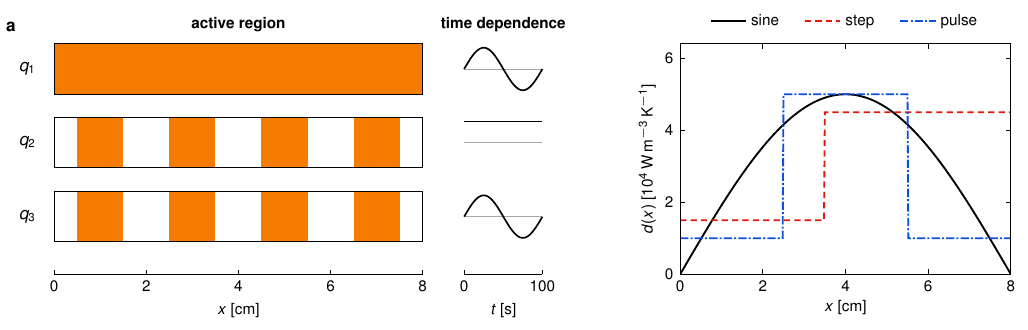}
\caption{Test problem conditions:
(a) active regions of sources $q_1$, $q_2$,
and $q_3$ and their respective temporal dependencies;
(b) sinusoidal, step, and pulse dissipation profiles.
The colored regions indicate
the action of the source, including during the
heat extraction half-cycle of $q_1$ and $q_3$.}
\label{fig:setup}
\end{figure*}

Figure~\ref{fig:setup} presents the
heating conditions and the dissipation
coefficient profiles.

\subsection{Inverse problem details and evaluation metrics}
\label{subsec:detalhes_inverso}

The synthetic temperatures are obtained using
the finite volume method, with $800$ cells and
Crank--Nicolson integration with a time step
of $0.05\,\mathrm{s}$.
The fields are stored at intervals
of $0.25\,\mathrm{s}$ and interpolated
at the sensor positions.
Using an independent solver
avoids generating the measurements with the same
discrete model employed in the inversion.
Gaussian noise with standard deviations
$\sigma=0.05$, $0.10$, and $0.20\,^\circ\mathrm C$
is added to the reference fields.
The reference case uses $N_x=81$ sensors
and $N_t=101$ observation times, with
a measurement interval of $1\,\mathrm{s}$.

The sensors are positioned at the midpoints
of the spatial intervals,

\begin{equation}
\ColumnMath[2.8em]{%
x_s\MathBreak
=
\left(s-\frac12\right)\frac{L}{N_x},\MathBreak
\qquad s=1,\ldots,N_x,%
}
\label{eq:posicoes_sensores}
\end{equation}

\noindent and the observation times include
the endpoints of the temporal interval.
The transformation employs the midpoint
quadrature rule, with $\omega_s=L/N_x$.
The eigenvalue problem uses an auxiliary basis of
$160$ cosine functions, with $60$ modes calculated.
For the comparison of the selection strategies in
Section~\ref{subsec:comparacao_truncamento} and for the
spectral and truncation-error verifications reported below,
$70$ modes are calculated with the same auxiliary basis; the
first $60$ are identical in the two calculations.
The integrals of the eigenvalue problem and the
matrices $\mathbf P^{(k)}$ are evaluated using
composite Gauss--Legendre quadrature with $4000$ points.

The properties in Eqs.~\eqref{eq:propriedades_teste}
and~\eqref{eq:transicao_propriedades} and the three profiles
of $d(x)$ satisfy assumptions (A1) and (A2), and $k(x)$ and
$w(x)$ are infinitely differentiable, so that all parts of
Lemma~\ref{lem:eig} apply. For this problem,
$\Phi=117.8\,\mathrm{s^{1/2}}$ and $e^{\mathcal V}=4.82$ in
Eq.~\eqref{eq:phi_v}, and the computed eigenvalues satisfy
$|\mu_i\Phi/[(i-1)\pi]-1|<3\times10^{-4}$ for
$10\le i\le70$, in agreement with Lemma~\ref{lem:eig}(c).

In each minimization, $N$ coefficients are estimated.
The initial distribution
$\hat d^{(0)}(x)=0.005\,\mathrm{s^{-1}}\,w(x)$
corresponds to the coefficients

\begin{equation}
\ColumnMath[2.8em]{%
\hat a_1^{(0)}\MathBreak
=
0.005\,\mathrm{s^{-1}}\sqrt W,\MathBreak
\qquad
\hat a_k^{(0)}=0,\MathBreak
\qquad k=2,\ldots,N.%
}
\label{eq:estimativa_inicial_coeficientes}
\end{equation}

The Levenberg--Marquardt damping parameter is initialized
at $\lambda_{\mathrm{LM},0}=10^{-3}$ and is reduced by a factor of ten when a step
is accepted and increased by the same factor
when it is rejected.
Minimization is terminated when
$|\Delta S|/S<10^{-10}$ over two iterations or
$\|\mathbf J^{\mathsf T}\mathbf r\|/
\|\mathbf J_0^{\mathsf T}\mathbf r_0\|<10^{-9}$,
with a maximum of $100$ iterations.

Reconstruction quality is evaluated
using the normalized root mean square error,

\begin{equation}
\ColumnMath[2.8em]{%
\operatorname{NRMSE}\MathBreak
=
\frac{
\left[
\dfrac{1}{N_e}
\sum_{r=1}^{N_e}
\left(
\hat d_N(x_r)-d(x_r)
\right)^2
\right]^{1/2}
}{
\max_r d(x_r)-\min_r d(x_r)
},%
}
\label{eq:nrmse}
\end{equation}

\noindent where $d(x_r)$ and $\hat d_N(x_r)$
are the exact and estimated values of the
dissipation coefficient.
The error is calculated at $N_e=401$ uniformly
spaced points, and its presentation as a percentage
corresponds to $100\operatorname{NRMSE}$.

\subsection{Verification of the direct solution}
\label{subsec:verificacao_direto}

The direct solution is verified using the exact
sinusoidal coefficient from Eq.~\eqref{eq:d_seno},
with the coupling matrix calculated from
Eq.~\eqref{eq:acoplamento_fonte}.
The comparison uses $81$ sensors and $101$ observation
times, with orders $N=10$, $20$, $30$, and $40$.
Crank--Nicolson time integration \cite{Crank1996, Ozisik2017} of the modal
system employs a time step of $0.25\,\mathrm{s}$.
Table~\ref{tab:verification} presents
the maximum absolute temperature difference
at the sensor positions and observation times.
Relative to the finite volume solution, this difference
decreases as $N$ increases, reaching
$4.8\times10^{-4}\,^\circ\mathrm C$
for $q_1$ and $0.064\,^\circ\mathrm C$
for $q_2$ and $q_3$, with $N=40$.
The slower convergence for the zonal sources
is consistent with their spatial discontinuities.

\begin{table*}[!tp]
\centering
\small
\caption{Verification of the direct solution:
maximum absolute temperature difference
[$^\circ\mathrm C$] relative to the finite
volume solution. The last column compares the
analytical solution with Crank--Nicolson, using
a time step of $0.25\,\mathrm{s}$, for the
orders examined.}
\label{tab:verification}
\setlength{\tabcolsep}{5pt}
\renewcommand{\arraystretch}{1.15}
\begin{adjustbox}{max width=\linewidth}
\begin{tabular}{lccccc}
\toprule
&
\multicolumn{4}{c}{\textbf{Modal vs.\ finite volume}}
& \textbf{Analytical vs.}
\\
\cmidrule(lr){2-5}
\textbf{Source}
& $\boldsymbol{N=10}$
& $\boldsymbol{N=20}$
& $\boldsymbol{N=30}$
& $\boldsymbol{N=40}$
& \textbf{Crank--Nicolson}
\\
\midrule
$q_1$
& $5.7\times10^{-2}$
& $6.0\times10^{-3}$
& $1.5\times10^{-3}$
& $4.8\times10^{-4}$
& $1.4\times10^{-4}$
\\
$q_2$
& $3.0$
& $2.4\times10^{-1}$
& $1.2\times10^{-1}$
& $6.4\times10^{-2}$
& $1.4\times10^{-4}$
\\
$q_3$
& $2.4$
& $2.4\times10^{-1}$
& $1.2\times10^{-1}$
& $6.4\times10^{-2}$
& $2.1\times10^{-4}$
\\
\bottomrule
\end{tabular}
\end{adjustbox}
\end{table*}

For a given modal order, the differences
between the analytical solution and Crank--Nicolson
do not exceed
$2.1\times10^{-4}\,^\circ\mathrm C$,
indicating that the discrepancies relative
to the finite volume solution arise predominantly
from the spatial representation and the
discretizations being compared. For the zonal sources,
the maximum difference with $N=40$ is
comparable to the lowest noise standard deviation
considered.

\subsection{A priori truncation order selection}
\label{subsec:resultados_truncamento_previo}

\begin{figure*}[!htbp]
\centering
\includegraphics[width=\textwidth]{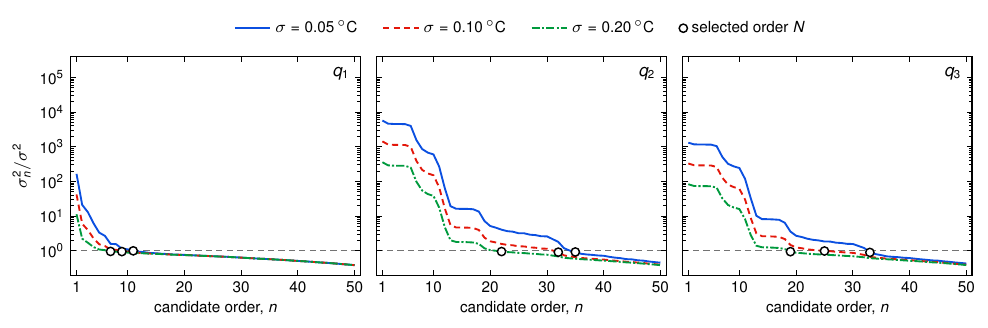}
\caption{A priori truncation selection for
the step coefficient: ratio $\sigma_n^2/\sigma^2$
as a function of the candidate order, for sources
$q_1$, $q_2$, and $q_3$, with $\sigma=0.05$,
$0.10$, and $0.20\,^\circ\mathrm C$.
The circles indicate the first order satisfying
the discrepancy condition; the dashed line corresponds
to the unit level. Reference case with
$81\times101$ measurements, one realization per
noise level.}
\label{fig:discrepancy}
\end{figure*}

In the reference case, $N_{\max}=50$ is adopted,
and the order is selected a priori using
Eq.~\eqref{eq:selecao_n}.
A total of $81\times101$ measurements are used, with
$\sigma=0.05$, $0.10$, and $0.20\,^\circ\mathrm C$,
considering one realization per noise level in
Figure~\ref{fig:discrepancy}.

Figure~\ref{fig:discrepancy} presents the ratio
$\sigma_n^2/\sigma^2$ for the step coefficient,
with $\sigma_n^2$ defined in
Eq.~\eqref{eq:discrepancia_previa}.
The circles indicate the first order satisfying
Eq.~\eqref{eq:selecao_n}, establishing the
dimension of the inverse problem before minimization.

Considering all profiles and realizations
of the reference case, source $q_1$ leads
to orders between $6$ and $15$, whereas $q_2$
and $q_3$ lead to orders between $19$ and $35$.
For the step profile and $\sigma=0.10\,^\circ\mathrm C$,
the selected orders are $N=9$, $31$--$32$,
and $24$--$26$, respectively.
The zonal sources produce temperature fields
with greater spatial content, requiring more
modes to represent the measurements at the
prescribed discrepancy level.
Since the same order is used in the expansion
of $\hat d_N(x)$, this behavior also
increases the number of estimated coefficients.

This behavior is described by
Proposition~\ref{prop:apriori}. For the sensor positions and
quadrature weights adopted, $\kappa_n$ in
Eq.~\eqref{eq:kappa} differs from $1-n/N_x$ by less than
$5\times10^{-5}$ for $n\le50$, so that each retained mode
removes $1/N_x=1.23\%$ of the noise variance from the
residual. The first order satisfying
Eq.~\eqref{eq:cruzamento_previo}, evaluated with the
noise-free reference fields, falls within the range of
orders selected in the ten realizations for all $27$
combinations of source, profile, and noise level. The
selected order is therefore determined by the decay of
$\beta_n^2$, that is, by the spatial content of the
temperature field.

Least-squares fits of $\log\beta_n^2$ against $\log n$ give
the decay exponent of Remark~\ref{rem:scaling} as
$p\approx2.9$ for the step and pulse profiles with $q_1$
($6\le n\le40$) and $p\approx3.0$ for the zonal sources
with the three profiles ($10\le n\le60$). For the sinusoidal
profile with $q_1$, in which neither the source nor the
coefficient is discontinuous, $p\approx4$. Since the
exponents obtained with $q_1$ and with the zonal sources are
similar for the discontinuous profiles, the higher orders
selected with $q_2$ and $q_3$ are associated with the larger
amplitude of the high-order components of the temperature
field, and not with a slower decay.

The highest selected order is $35$, below
the limit $N_{\max}=50$.
Therefore, this limit does not determine the
selections in the reference case.
Reducing $\sigma$ allows lower-amplitude
components to be retained and generally increases $N$.
According to Eq.~\eqref{eq:escala_ordem},
$N\propto\sigma^{-1/p}$. In Table~\ref{tab:estimates},
reducing $\sigma$ by a factor of four increases the mean
order from $7.0$ to $11.1$ for the step with $q_1$ and from
$21.0$ to $35.0$ with $q_2$, which corresponds to exponents
of $0.33$ and $0.37$, close to $1/p\simeq1/3$. For the
sinusoidal profile with $q_1$, the mean order increases only
from $6.9$ to $9.0$, an exponent of $0.19$, in line with the
larger value $p\approx4$.

\subsection{Reconstruction of the dissipation coefficient}
\label{subsec:resultados_reconstrucao}

\begin{figure*}[!htbp]
\centering
\includegraphics[width=\textwidth]{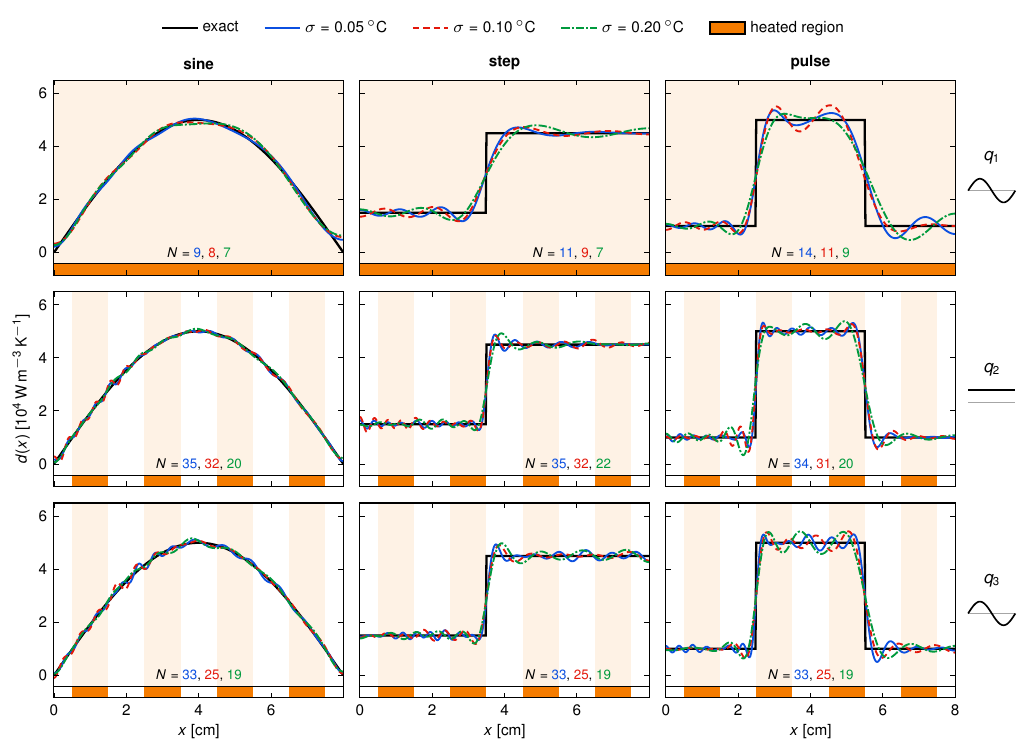}
\caption{Reconstructions $\hat d_N(x)$ for sources
$q_1$, $q_2$, and $q_3$ in the rows and the sinusoidal,
step, and pulse profiles in the columns. Results
are presented for $\sigma=0.05$, $0.10$, and
$0.20\,^\circ\mathrm C$, with $N$ selected a priori. The black curve represents $d(x)$.
The curves correspond to the realization with
random seed equal to one; the statistics are presented
in Table~\ref{tab:estimates}.}
\label{fig:estimates}
\end{figure*}

The reconstructions use $81\times101$ measurements,
with $\sigma=0.05$, $0.10$, and $0.20\,^\circ\mathrm C$.
The statistics are calculated over ten
noise realizations, reporting the mean and
population standard deviation of the NRMSE.
The reconstruction curves use the
realization with random seed equal to one.
Figure~\ref{fig:estimates} presents the functions
$\hat d_N(x)$ obtained using a priori selection,
and Table~\ref{tab:estimates} summarizes the statistics
of the ten realizations.
The sinusoidal profile exhibits mean errors between
$1.11\%$ and $2.84\%$.
For the step and pulse profiles, the ranges
are $5.63\%$ to $12.04\%$ and $8.09\%$
to $15.60\%$, respectively.

\begin{table*}[!tp]
\centering
\small
\caption{Reconstruction of $d(x)$ using a priori
selection: range of selected orders,
mean NRMSE and population standard deviation over
ten realizations, and mean number of
Levenberg--Marquardt iterations (It.).}
\label{tab:estimates}
\setlength{\tabcolsep}{5pt}
\renewcommand{\arraystretch}{1.15}
\begin{adjustbox}{max width=\linewidth}
\begin{tabular}{llccccccccc}
\toprule
& &
\multicolumn{3}{c}{$\sigma=0.05\,^\circ\mathrm C$}
&
\multicolumn{3}{c}{$\sigma=0.10\,^\circ\mathrm C$}
&
\multicolumn{3}{c}{$\sigma=0.20\,^\circ\mathrm C$}
\\
\cmidrule(lr){3-5}
\cmidrule(lr){6-8}
\cmidrule(lr){9-11}
\textbf{Source}
& $\boldsymbol{d(x)}$
& $\boldsymbol{N}$ & \textbf{NRMSE [\%]} & \textbf{It.}
& $\boldsymbol{N}$ & \textbf{NRMSE [\%]} & \textbf{It.}
& $\boldsymbol{N}$ & \textbf{NRMSE [\%]} & \textbf{It.}
\\
\midrule
$q_1$ & sinusoidal
& 9 & 1.77 $\pm$ 0.01 & 5.0
& 7--8 & 2.56 $\pm$ 0.10 & 5.0
& 6--7 & 2.84 $\pm$ 0.40 & 5.5
\\
& step
& 11--12 & 9.66 $\pm$ 0.11 & 5.0
& 9 & 10.56 $\pm$ 0.01 & 5.0
& 7 & 12.04 $\pm$ 0.02 & 5.2
\\
& pulse
& 14--15 & 12.83 $\pm$ 0.15 & 5.0
& 11 & 13.88 $\pm$ 0.01 & 5.1
& 9 & 15.60 $\pm$ 0.02 & 6.0
\\
\midrule
$q_2$ & sinusoidal
& 35 & 1.11 $\pm$ 0.24 & 5.0
& 31--32 & 1.41 $\pm$ 0.37 & 5.0
& 20 & 1.17 $\pm$ 0.07 & 5.0
\\
& step
& 35 & 5.63 $\pm$ 0.10 & 5.0
& 31--32 & 6.19 $\pm$ 0.18 & 5.0
& 20--22 & 7.34 $\pm$ 0.15 & 5.0
\\
& pulse
& 34--35 & 8.09 $\pm$ 0.03 & 5.0
& 31 & 8.55 $\pm$ 0.09 & 5.0
& 20 & 10.99 $\pm$ 0.03 & 5.0
\\
\midrule
$q_3$ & sinusoidal
& 33 & 1.84 $\pm$ 0.38 & 5.0
& 24--25 & 1.82 $\pm$ 0.24 & 5.3
& 19 & 1.83 $\pm$ 0.36 & 6.0
\\
& step
& 33 & 6.34 $\pm$ 0.27 & 5.0
& 24--26 & 7.06 $\pm$ 0.17 & 5.1
& 19 & 8.02 $\pm$ 0.31 & 6.0
\\
& pulse
& 33 & 8.35 $\pm$ 0.11 & 5.0
& 24--25 & 9.92 $\pm$ 0.08 & 5.9
& 19 & 11.34 $\pm$ 0.07 & 6.0
\\
\bottomrule
\end{tabular}
\end{adjustbox}
\end{table*}

For the smooth profile, the estimates follow
the exact distribution along the plate,
with small differences near the boundaries
and in the curvature.
These differences are explained by
Proposition~\ref{prop:neumann}: for the profile of
Eq.~\eqref{eq:d_seno}, the ratio $d/w$ has a non-zero
derivative at $x=0$ and $x=L$, whereas all eigenfunctions
have zero slope at the ends, so that the coefficients decay
only as $k^{-2}$ and the truncation error as $n^{-3/2}$.
For the properties of the test problem,
$\operatorname{NRMSE}(d_n)\,n^{3/2}$ lies between $0.48$ and
$0.57$ for $n=10,\ldots,70$, against
$2/(\sqrt3\,\pi)\approx0.37$ in
Eq.~\eqref{eq:taxa_seno}.
For the discontinuous profiles, the locations
of the transitions are recovered, but
oscillations and local smoothing of the jumps
appear, consistent with representation by a
truncated expansion in smooth functions.
The pulse exhibits a higher normalized error
than the step under the conditions compared,
involving the reconstruction of two transitions
and the plateau between them.
Theorem~\ref{thm:step} quantifies this behavior: the
truncation error decays as $C_d/\sqrt N$, with
$C_d=1/\pi\approx31.8\%$ for the step and
$C_d=\sqrt2/\pi\approx45.0\%$ for the pulse in a uniform
medium. For the heterogeneous properties considered,
Eq.~\eqref{eq:Cd_het} gives $31.4\%$ and $46.6\%$, and the
direct evaluation of $\operatorname{NRMSE}(d_n)\sqrt n$ for
$n=30,\ldots,70$ gives $32.0\%$ (range $31.5$--$32.4$) and
$47.4\%$ (range $46.5$--$48.2$). The ratio between these
constants accounts for the higher error of the pulse at
comparable orders.

The influence of the source is more evident
for the discontinuous profiles.
For $\sigma=0.10\,^\circ\mathrm C$,
the mean error for the step decreases from $10.56\%$
with $q_1$ to $6.19\%$ with $q_2$
and $7.06\%$ with $q_3$.
For the pulse, the corresponding values
are $13.88\%$, $8.55\%$, and $9.92\%$.
The higher orders selected for the zonal sources
allow narrower transitions to be represented.
The difference between $q_2$ and $q_3$ also involves
the temporal evolution of heating, since
these sources have the same spatial distribution.

Table~\ref{tab:floor} compares the truncation error of
$d_N$, obtained with the exact coefficients of
Eq.~\eqref{eq:coeficientes_d} at the orders selected
a priori, with the mean errors of
Table~\ref{tab:estimates}. By
Proposition~\ref{prop:decomp}, the truncation term cannot be
reduced by any estimate of order $N$. With $q_1$, the
reported errors coincide with the truncation error to within
about $0.1$ percentage point in all cases: these
reconstructions are limited by the truncation order, and not
by the measurement errors or by the estimator. With the
zonal sources and the discontinuous profiles, the excess is
at most $0.8$ percentage point. Only the sinusoidal profile
with $q_2$ and $q_3$ exhibits errors dominated by the
estimation term. The lower errors obtained with the zonal
sources for the step and pulse profiles therefore result
from the higher orders selected, through the factor
$C_d/\sqrt N$ of Theorem~\ref{thm:step}.

\begin{table*}[!tp]
\centering
\small
\caption{Truncation error
$100\operatorname{NRMSE}(d_N)$ at the orders selected
a priori (first value; a range when the selected order
varies among the ten realizations) and mean NRMSE of
Table~\ref{tab:estimates} (second value), in \%.}
\label{tab:floor}
\setlength{\tabcolsep}{7pt}
\renewcommand{\arraystretch}{1.15}
\begin{adjustbox}{max width=\linewidth}

\begin{tabular}{llccc}
\toprule
\textbf{Source}
& $\boldsymbol{d(x)}$
& $\sigma=0.05\,^\circ\mathrm C$
& $\sigma=0.10\,^\circ\mathrm C$
& $\sigma=0.20\,^\circ\mathrm C$
\\
\midrule
$q_1$ & sinusoidal
& 1.75 / 1.77 & 2.44--2.66 / 2.56 & 2.66--3.97 / 2.84
\\
& step
& 9.27--9.63 / 9.66 & 10.51 / 10.56 & 11.99 / 12.04
\\
& pulse
& 12.53--12.86 / 12.83 & 13.82 / 13.88 & 15.49 / 15.60
\\
\midrule
$q_2$ & sinusoidal
& 0.25 / 1.11 & 0.29 / 1.41 & 0.57 / 1.17
\\
& step
& 5.36 / 5.63 & 5.58--5.76 / 6.19 & 6.86--7.14 / 7.34
\\
& pulse
& 8.02--8.03 / 8.09 & 8.40 / 8.55 & 10.88 / 10.99
\\
\midrule
$q_3$ & sinusoidal
& 0.27 / 1.84 & 0.39--0.44 / 1.82 & 0.59 / 1.83
\\
& step
& 5.56 / 6.34 & 6.26--6.56 / 7.06 & 7.33 / 8.02
\\
& pulse
& 8.09 / 8.35 & 9.68--9.82 / 9.92 & 11.12 / 11.34
\\
\bottomrule
\end{tabular}
\end{adjustbox}
\end{table*}

For the step and pulse profiles, reducing the noise
monotonically decreases the mean error for
all three sources.
For the step with $q_2$, the error decreases from
$7.34\%$ to $5.63\%$ when $\sigma$
is reduced from $0.20$ to
$0.05\,^\circ\mathrm C$.
For the sinusoidal profile, the trend is not
strictly monotonic: with $q_2$, the errors are
$1.11\%$, $1.41\%$, and $1.17\%$
at increasing noise levels;
with $q_3$, they remain close to $1.8\%$.
Thus, the selected order provides adequate
reconstructions of the smooth profile in the
cases examined.
In these cases, the truncation error at the selected
orders, between $0.25\%$ and $0.59\%$
(Table~\ref{tab:floor}), is well below the reconstruction
error, so that, by Proposition~\ref{prop:decomp} and
Remark~\ref{rem:noise}, the additional coefficients
contribute only through the propagation of the measurement
errors.

Fitting in the transformed domain is also
monitored through the ratio $S_B/(NN_t)$.
This quantity indicates the magnitude of the
residual relative to the modal noise scale.
For the sensors and weights adopted, the diagonal
elements of $\boldsymbol\Sigma/\bar\sigma^2$ lie between
$0.99$ and $1.02$ and the off-diagonal elements are below
$0.06$ in modulus, so that the common scale is close to the
full covariance of Eq.~\eqref{eq:covariancia_modal}
(Remark~\ref{rem:noise}).

The mean values of $S_B/(NN_t)$ range from $0.984$
to $1.025$ for the configurations in
Table~\ref{tab:estimates}.
The fitted residuals remain close
to the modal noise scale, even in cases
with larger errors in $\hat d_N(x)$.
This result shows that the reconstruction error
is governed by the spatial representation of the
unknown function, and not by the quality of the
temperature fit.
This is the content of Proposition~\ref{prop:decomp}:
the last term of Eq.~\eqref{eq:pitagoras} does not depend on
the estimated coefficients.

\subsection{Minimization behavior and comparison of estimators}
\label{subsec:comparacao_estimadores}

Levenberg--Marquardt completes the estimation
in five or six iterations in the $270$ runs
reported in Table~\ref{tab:estimates}.
Gauss--Newton was evaluated for one realization
of each combination of source and profile,
with $\sigma=0.10\,^\circ\mathrm C$,
yielding the same number of iterations
in the corresponding comparisons.

\begin{figure*}[!htbp]
\centering
\includegraphics[width=\textwidth]{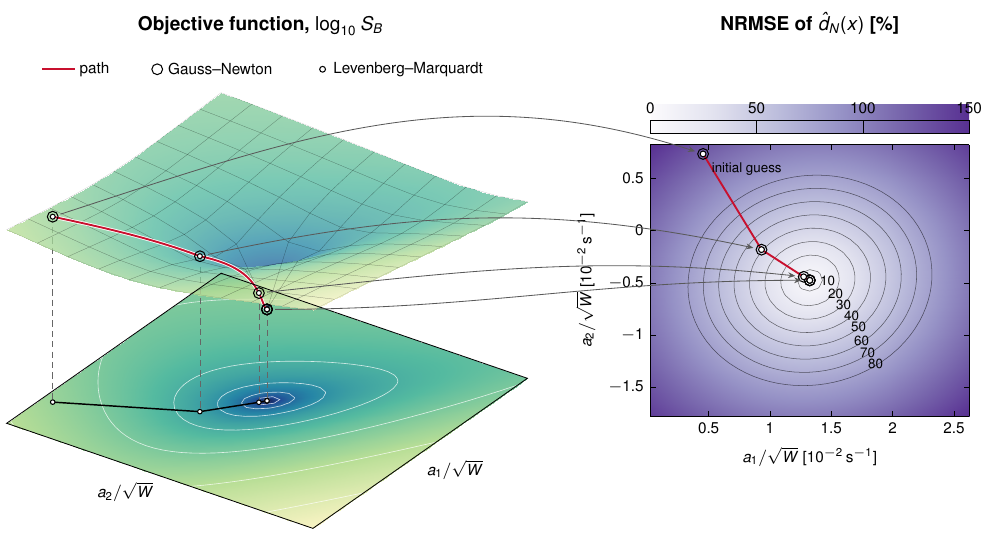}
\caption{Reduced problem in the
$(\hat a_1,\hat a_2)$ plane for the step coefficient,
source $q_2$, and $\sigma=0.10\,^\circ\mathrm C$,
with the remaining coefficients fixed at their
projected values of $d(x)$. On the left, the rescaled
surface of $\log_{10}S_B$ and minimization trajectories;
on the right, the NRMSE in the same plane.
The larger and smaller circles represent
the Gauss--Newton and Levenberg--Marquardt
iterations, respectively. The arrows connect
the corresponding positions in the panels.}
\label{fig:surface}
\end{figure*}

Figure~\ref{fig:surface} presents a
reduced problem for the step profile with $q_2$,
$81\times101$ measurements, and $\sigma=0.10\,^\circ\mathrm C$,
in which only $\hat a_1$ and $\hat a_2$ vary,
while the remaining coefficients are kept
at the values obtained from the projection of $d(x)$.
The functional $\log_{10}S_B$, rescaled for
visualization, the NRMSE, and the trajectories
of the two estimators are presented.
These trajectories correspond to the optimization
restricted to two parameters.

Both methods approach the region of
lower residual, which is also associated
with a lower reconstruction error.
The final NRMSE of the reduced problem is
approximately $5.58\%$, whereas the
full inversion of the same realization
yields $6.34\%$.
The difference is consistent with fixing
the remaining coefficients at their projected values.
The error of the reduced problem coincides, to the
reported precision, with the truncation error of order
$N=32$, the lower limit of the corresponding range in
Table~\ref{tab:floor}, and the difference with respect to
the full inversion corresponds to the first term of
Eq.~\eqref{eq:pitagoras}.
Theorem~\ref{thm:ident} provides the corresponding
statement in the full parameter space: under conditions
(I1) and (I2), the latter being automatic in formulation B,
$\mathbf J^{\mathsf T}\mathbf J$ is positive definite and a
minimizer with zero residual is isolated
(Remark~\ref{rem:ident}).

The computational cost of the estimators is evaluated
using elapsed time and the number of model calls.
The Levenberg--Marquardt calls include
evaluations with the Jacobian.

Table~\ref{tab:nminimize} presents a
comparison with \texttt{NMinimize} \cite{WolframNMinimize}, using
Nelder--Mead and a model implementation
in Wolfram Language \cite{Mathematica15}. The comparison uses
the step coefficient, $81\times101$ measurements,
$\sigma=0.10\,^\circ\mathrm C$, and the same
data realization for both estimators.
All $N$ coefficients are estimated,
with the following constraints in the
\texttt{NMinimize} execution:

\begin{equation}
\ColumnMath[2.8em]{%
-0.08\,\mathrm{s^{-1}}\MathBreak
\leq
\frac{\hat a_k}{\sqrt W}\MathBreak
\leq
0.08\,\mathrm{s^{-1}},\MathBreak
\qquad k=1,\ldots,N.%
}
\label{eq:restricoes_nminimize}
\end{equation}

The final functional values agree
with those obtained by Levenberg--Marquardt,
with relative differences on the order
of $10^{-14}$ to $10^{-11}$,
and the reconstruction errors agree
to the reported precision.

The differences between the coefficients obtained
by the estimators are evaluated using
$\max_k|\Delta\hat a_k|/\sqrt W$, in
$\mathrm{s^{-1}}$, where
$\Delta\hat a_k$ is the difference between their
final estimates.

\begin{table*}[!tp]
\centering
\small
\caption{Comparison between Levenberg--Marquardt
(LM), with sensitivities obtained from the spectral solution,
and \texttt{NMinimize} with Nelder--Mead:
model calls, elapsed time,
relative difference between the minima, and
maximum difference between coefficients divided
by $\sqrt W$, in $\mathrm{s^{-1}}$.
Step coefficient,
$\sigma=0.10\,^\circ\mathrm C$,
same data realization, and constraints
from Eq.~\eqref{eq:restricoes_nminimize}
in \texttt{NMinimize}.
LM calls include evaluations
with the Jacobian.}
\label{tab:nminimize}
\setlength{\tabcolsep}{5pt}
\renewcommand{\arraystretch}{1.15}
\begin{adjustbox}{max width=\linewidth}
\begin{tabular}{lclrrccc}
\toprule
\textbf{Source}
& $\boldsymbol{N}$
& \textbf{Estimator}
& \textbf{Calls}
& \textbf{Time [s]}
& $|\Delta S|/S$
& $\max_k|\Delta\hat a_k|/\sqrt W$
& \textbf{NRMSE [\%]}
\\
\midrule
$q_1$ & 9
& LM (spectral sensitivities)
& 11 & 0.018 & -- & -- & 10.56
\\
& & \texttt{NMinimize} (Nelder--Mead)
& 5327 & 86
& $6.7\times10^{-14}$
& $9.6\times10^{-13}$
& 10.56
\\
\midrule
$q_2$ & 32
& LM (spectral sensitivities)
& 11 & 0.207 & -- & -- & 6.34
\\
& & \texttt{NMinimize} (Nelder--Mead)
& 26193 & 1196
& $8.4\times10^{-12}$
& $3.8\times10^{-9}$
& 6.34
\\
\midrule
$q_3$ & 25
& LM (spectral sensitivities)
& 11 & 0.095 & -- & -- & 7.04
\\
& & \texttt{NMinimize} (Nelder--Mead)
& 24034 & 980
& $1.1\times10^{-12}$
& $1.2\times10^{-9}$
& 7.04
\\
\bottomrule
\end{tabular}
\end{adjustbox}
\end{table*}

Levenberg--Marquardt requires $11$ model calls,
including evaluations of the residuals
and the Jacobian, whereas Nelder--Mead requires
$5327$, $26193$, and $24034$ evaluations for
$q_1$, $q_2$, and $q_3$.
The elapsed times range from $0.018$ to
$0.207\,\mathrm{s}$ for LM and from $86$
to $1196\,\mathrm{s}$ for \texttt{NMinimize}.

The agreement between the two estimators reinforces
the consistency of the estimates.
Since the computational environments and implementations
are different, the number of model evaluations provides
the most direct comparison, and its reduction highlights
the benefit of using sensitivities
during minimization.

\subsection{Influence of measurement resolution and source amplitude}
\label{subsec:resolucao_amplitude}

\begin{figure*}[!htbp]
\centering
\includegraphics[width=\textwidth]{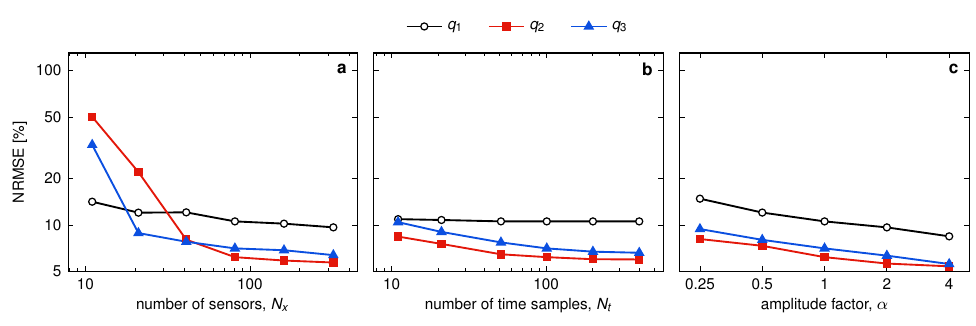}
\caption{Influence of measurements and source amplitude
on the mean NRMSE of the step coefficient reconstruction,
with $\sigma=0.10\,^\circ\mathrm C$:
(a) number of sensors, with $N_t=101$;
(b) number of observation times, with $N_x=81$;
(c) amplitude factor, with $81\times101$ measurements.
Means over ten realizations, with a priori
truncation selection in each configuration.}
\label{fig:data}
\end{figure*}

Figure~\ref{fig:data} presents the influence
of the number of measurements and source amplitude
on the reconstruction of the step coefficient,
keeping $\sigma=0.10\,^\circ\mathrm C$.
Tables~\ref{tab:sensors}--\ref{tab:amplitude}
summarize the results for the three dissipation
profiles. In each configuration, the order $N$
is selected again a priori using
Eq.~\eqref{eq:selecao_n}.
Thus, variations in the error also include
changes in the number of estimated coefficients.

In the spatial resolution study, $N_t=101$ is kept fixed;
in the temporal resolution study, $N_x=81$; and in the
amplitude study, $81\times101$ measurements are used.
In the tables, the values in
parentheses correspond to the mean selected
order, rounded to an integer.

In the analyses with different numbers
of sensors, the maximum candidate order is
limited by

\begin{equation}
\ColumnMath[2.8em]{%
N_{\max}\MathBreak
=
\min(50,N_x-1).%
}
\label{eq:limite_truncamento_numerico}
\end{equation}

In the amplitude study, each source is
multiplied by the dimensionless factor $\alpha$,
while keeping the noise standard deviation fixed.
For a given $d(x)$, the linearity of the direct
problem with respect to the source causes an increase
in $\alpha$ to increase the thermal response
and the signal-to-noise ratio.

\begin{table*}[!tp]
\centering
\small
\caption{Influence of the number of spatial sensors:
mean NRMSE [\%] over ten realizations, with $N_t=101$
and $\sigma=0.10\,^\circ\mathrm C$.
In parentheses, the mean order selected
a priori, rounded to an integer.}
\label{tab:sensors}
\setlength{\tabcolsep}{5pt}
\renewcommand{\arraystretch}{1.15}
\begin{adjustbox}{max width=\linewidth}
\begin{tabular}{rccccccccc}
\toprule
& \multicolumn{3}{c}{$q_1$}
& \multicolumn{3}{c}{$q_2$}
& \multicolumn{3}{c}{$q_3$}
\\
\cmidrule(lr){2-4}
\cmidrule(lr){5-7}
\cmidrule(lr){8-10}
$\boldsymbol{N_x}$
& \textbf{Sinusoidal} & \textbf{Step} & \textbf{Pulse}
& \textbf{Sinusoidal} & \textbf{Step} & \textbf{Pulse}
& \textbf{Sinusoidal} & \textbf{Step} & \textbf{Pulse}
\\
\midrule
11
& 4.2 (6) & 14.2 (6) & 17.6 (8)
& 41.1 (10) & 50.2 (10) & 48.1 (10)
& 25.7 (10) & 33.0 (10) & 32.7 (10)
\\
21
& 2.7 (7) & 12.0 (7) & 15.6 (9)
& 13.3 (16) & 22.0 (17) & 18.7 (15)
& 2.6 (15) & 8.9 (15) & 12.8 (15)
\\
41
& 2.7 (7) & 12.1 (8) & 15.4 (9)
& 3.3 (25) & 8.0 (25) & 10.3 (24)
& 2.0 (20) & 7.8 (21) & 11.0 (20)
\\
81
& 2.6 (8) & 10.6 (9) & 13.9 (11)
& 1.4 (31) & 6.2 (32) & 8.5 (31)
& 1.8 (25) & 7.1 (25) & 9.9 (25)
\\
161
& 1.8 (9) & 10.2 (10) & 13.4 (13)
& 1.1 (33) & 5.9 (33) & 8.2 (33)
& 2.1 (31) & 6.9 (31) & 8.7 (31)
\\
321
& 1.8 (9) & 9.7 (11) & 12.7 (14)
& 0.9 (34) & 5.7 (34) & 8.1 (34)
& 1.8 (33) & 6.4 (33) & 8.3 (33)
\\
\bottomrule
\end{tabular}
\end{adjustbox}
\end{table*}

Increasing spatial resolution is particularly
relevant for $q_2$. The mean error for the step
decreases from $50.2\%$, with $N_x=11$, to
$22.0\%$, with $N_x=21$, reaching $8.0\%$
and $6.2\%$ with $41$ and $81$ sensors.
The improvement in the reconstruction of $d(x)$ depends on the sensitivity of
the temperatures to the additional coefficients
and the accuracy of the measurement transformation.

The advantage of zonal excitation therefore
depends on sampling capable of representing
the resulting temperature field. With few sensors,
quadrature accuracy and the constraint
$N\leq N_x-1$ restrict the modal information
available for estimating the coefficients
$\hat a_k$.
By Proposition~\ref{prop:apriori}, each retained mode
removes a fraction $1/N_x$ of the noise variance from the
residual, and Remark~\ref{rem:ident} requires $N\le N_x$ for
the sampling matrix to have full rank. With $N_x=11$ and the
zonal sources, the mean selected order coincides with the
limit of Eq.~\eqref{eq:limite_truncamento_numerico}, so that
the order is set by the sampling and not by the decay of
$\beta_n^2$ in Eq.~\eqref{eq:cruzamento_previo}.

This behavior depends on the source and profile.
With $21$ sensors, $q_3$ already provides a lower
error for the step than $q_1$, whereas $q_2$
still exhibits a higher error.
With $41$ sensors, the zonal sources outperform
$q_1$ for the step and pulse profiles, but $q_2$
still exhibits a higher error for the sinusoidal profile:
$3.3\%$, compared with $2.7\%$ for $q_1$.
Beyond $81$ sensors, the improvements become
smaller. For the step with $q_2$, the error
decreases from $6.2\%$ to $5.7\%$ when
$N_x$ increases to $321$.

\begin{table*}[!tp]
\centering
\small
\caption{Influence of the number of observation times:
mean NRMSE [\%] over ten realizations, with $N_x=81$
and $\sigma=0.10\,^\circ\mathrm C$.
In parentheses, the mean order selected
a priori, rounded to an integer.}
\label{tab:samples}
\setlength{\tabcolsep}{5pt}
\renewcommand{\arraystretch}{1.15}
\begin{adjustbox}{max width=\linewidth}
\begin{tabular}{rccccccccc}
\toprule
& \multicolumn{3}{c}{$q_1$}
& \multicolumn{3}{c}{$q_2$}
& \multicolumn{3}{c}{$q_3$}
\\
\cmidrule(lr){2-4}
\cmidrule(lr){5-7}
\cmidrule(lr){8-10}
$\boldsymbol{N_t}$
& \textbf{Sinusoidal} & \textbf{Step} & \textbf{Pulse}
& \textbf{Sinusoidal} & \textbf{Step} & \textbf{Pulse}
& \textbf{Sinusoidal} & \textbf{Step} & \textbf{Pulse}
\\
\midrule
11
& 2.6 (7) & 10.9 (9) & 14.3 (11)
& 3.9 (31) & 8.4 (31) & 9.6 (30)
& 5.4 (24) & 10.4 (24) & 11.5 (24)
\\
21
& 2.5 (8) & 10.8 (9) & 14.3 (11)
& 2.9 (31) & 7.5 (32) & 9.1 (31)
& 3.5 (24) & 9.0 (25) & 10.5 (24)
\\
51
& 2.5 (8) & 10.6 (9) & 14.0 (11)
& 1.9 (32) & 6.5 (32) & 8.7 (31)
& 2.6 (25) & 7.7 (25) & 10.2 (25)
\\
101
& 2.6 (8) & 10.6 (9) & 13.9 (11)
& 1.4 (31) & 6.2 (32) & 8.5 (31)
& 1.8 (25) & 7.1 (25) & 9.9 (25)
\\
201
& 2.6 (8) & 10.6 (9) & 13.9 (11)
& 1.1 (31) & 6.0 (32) & 8.5 (31)
& 1.3 (25) & 6.7 (25) & 9.8 (25)
\\
401
& 2.6 (7) & 10.6 (9) & 13.9 (11)
& 0.8 (31) & 6.0 (31) & 8.4 (31)
& 1.0 (25) & 6.6 (25) & 9.7 (25)
\\
\bottomrule
\end{tabular}
\end{adjustbox}
\end{table*}

Temporal refinement exhibits a different effect.
For $q_1$, the errors vary little between $11$
and $401$ observation times. For the sinusoidal profile with
zonal sources, the error decreases from $3.9\%$ to
$0.8\%$ with $q_2$ and from $5.4\%$ to
$1.0\%$ with $q_3$.
For the step, the corresponding reductions
are from $8.4\%$ to $6.0\%$ and from
$10.4\%$ to $6.6\%$.
The mean orders remain nearly constant throughout
this study, indicating that the improvement
arises mainly from the additional information
available for fitting the retained coefficients.
Increasing $N_t$, however, does not eliminate
the spatial limitation in representing jumps.
These trends follow from the results of
Section~\ref{sec:problema_inverso}. The margin in
Eq.~\eqref{eq:cruzamento_previo} does not depend on $N_t$,
which explains the nearly constant orders. By
Remark~\ref{rem:noise}, the contribution of the measurement
errors decreases as $N_t^{-1/2}$: for the sinusoidal profile
with the zonal sources, in which this contribution is
dominant, the errors decrease by factors of $4.8$ and $5.4$
between $N_t=11$ and $401$, against
$(401/11)^{1/2}\approx6.0$. With $q_1$, the errors are at
the truncation level of Table~\ref{tab:floor} and are
therefore insensitive to $N_t$. For the step, the errors
approach the values $C_d/\sqrt N$ of
Theorem~\ref{thm:step}, which depend only on the order.

\begin{table*}[!tp]
\centering
\small
\caption{Influence of the amplitude factor $\alpha$:
mean NRMSE [\%] over ten realizations, with
$81\times101$ measurements and
$\sigma=0.10\,^\circ\mathrm C$.
In parentheses, the mean order selected
a priori, rounded to an integer.}
\label{tab:amplitude}
\setlength{\tabcolsep}{5pt}
\renewcommand{\arraystretch}{1.15}
\begin{adjustbox}{max width=\linewidth}
\begin{tabular}{rccccccccc}
\toprule
& \multicolumn{3}{c}{$q_1$}
& \multicolumn{3}{c}{$q_2$}
& \multicolumn{3}{c}{$q_3$}
\\
\cmidrule(lr){2-4}
\cmidrule(lr){5-7}
\cmidrule(lr){8-10}
$\boldsymbol{\alpha}$
& \textbf{Sinusoidal} & \textbf{Step} & \textbf{Pulse}
& \textbf{Sinusoidal} & \textbf{Step} & \textbf{Pulse}
& \textbf{Sinusoidal} & \textbf{Step} & \textbf{Pulse}
\\
\midrule
0.25
& 4.5 (5) & 14.8 (5) & 19.6 (6)
& 1.7 (17) & 8.1 (18) & 11.6 (16)
& 2.0 (13) & 9.4 (13) & 13.3 (14)
\\
0.5
& 2.8 (7) & 12.0 (7) & 15.6 (9)
& 1.2 (20) & 7.3 (21) & 11.0 (20)
& 1.8 (19) & 8.0 (19) & 11.3 (19)
\\
1
& 2.6 (8) & 10.6 (9) & 13.9 (11)
& 1.4 (31) & 6.2 (32) & 8.5 (31)
& 1.8 (25) & 7.1 (25) & 9.9 (25)
\\
2
& 1.8 (9) & 9.7 (11) & 12.8 (14)
& 1.1 (35) & 5.6 (35) & 8.1 (34)
& 1.8 (33) & 6.3 (33) & 8.3 (33)
\\
4
& 1.3 (11) & 8.4 (14) & 11.3 (17)
& 1.6 (44) & 5.4 (45) & 7.2 (44)
& 1.3 (38) & 5.6 (38) & 7.8 (38)
\\
\bottomrule
\end{tabular}
\end{adjustbox}
\end{table*}

When $\alpha$ increases from $0.25$ to $4$,
the error for the step decreases from $14.8\%$ to
$8.4\%$ with $q_1$, from $8.1\%$ to
$5.4\%$ with $q_2$, and from $9.4\%$ to
$5.6\%$ with $q_3$.
The pulse exhibits a similar trend.
According to Eq.~\eqref{eq:escala_ordem},
$N\propto\alpha^{1/p}$: increasing $\alpha$ by a factor of
sixteen increases the mean order from $5$ to $14$ for the
step with $q_1$ and from $18$ to $45$ with $q_2$,
corresponding to exponents of $0.37$ and $0.33$, close to
$1/p\simeq1/3$.
For the sinusoidal profile with $q_2$, however, the error increases
from $1.1\%$ to $1.6\%$ between $\alpha=2$
and $4$, while the mean order increases from $35$
to $44$.
This result shows that including
additional coefficients does not ensure
greater accuracy for a smooth profile.
By Proposition~\ref{prop:neumann}, the truncation
error of the sinusoidal profile at these orders is already
well below the reconstruction error
(Table~\ref{tab:floor}), so that the additional coefficients
increase only the first term of
Eq.~\eqref{eq:pitagoras}.

At the highest amplitude, the mean orders for $q_2$
reach $44$--$45$, approaching the admissible limit.
The relationship between the number
of sensors and the selected order must
therefore be examined together with
quadrature accuracy and the characteristics
of the excitation.

\subsection{Computational cost of the temporal solution and functionals}
\label{subsec:custo_computacional}

The comparison of the temporal solution uses the
step coefficient, $81\times101$ measurements,
$\sigma=0.10\,^\circ\mathrm C$, and the reference
realization, with random seed equal to one.

Table~\ref{tab:solver} compares the analytical
solution in time with Crank--Nicolson in
Levenberg--Marquardt estimation.
For the step coefficient, both approaches require
five iterations and yield the same NRMSE values
to the reported precision.
The maximum difference between the estimated
coefficients, divided by $\sqrt W$, is less than
$5\times10^{-7}\,\mathrm{s^{-1}}$.

\begin{table*}[!tp]
\centering
\small
\caption{Comparison of the analytical solution in time
with Crank--Nicolson in Levenberg--Marquardt
estimation: iterations, elapsed time,
time ratio, NRMSE, and maximum difference
between estimated coefficients divided by
$\sqrt W$, in $\mathrm{s^{-1}}$.
Step coefficient, $81\times101$ measurements,
$\sigma=0.10\,^\circ\mathrm C$,
and reference realization.}
\label{tab:solver}
\setlength{\tabcolsep}{5pt}
\renewcommand{\arraystretch}{1.15}
\begin{adjustbox}{max width=\linewidth}
\begin{tabular}{lccccccccc}
\toprule
& &
\multicolumn{2}{c}{\textbf{Analytical}}
&
\multicolumn{2}{c}{\textbf{Crank--Nicolson}}
& &
\multicolumn{2}{c}{\textbf{NRMSE [\%]}}
&
\\
\cmidrule(lr){3-4}
\cmidrule(lr){5-6}
\cmidrule(lr){8-9}
\textbf{Source}
& $\boldsymbol{N}$
& \textbf{It.}
& \textbf{Time [ms]}
& \textbf{It.}
& \textbf{Time [ms]}
& \textbf{CN/analytical ratio}
& \textbf{Analytical}
& \textbf{CN}
& $\max_k|\Delta\hat a_k|/\sqrt W$
\\
\midrule
$q_1$ & 9
& 5 & 15.0
& 5 & 176.2
& 11.8 & 10.56 & 10.56
& $2.7\times10^{-7}$
\\
$q_2$ & 32
& 5 & 155.4
& 5 & 457.9
& 2.9 & 6.34 & 6.34
& $3.3\times10^{-8}$
\\
$q_3$ & 25
& 5 & 98.4
& 5 & 361.0
& 3.7 & 7.04 & 7.04
& $4.9\times10^{-7}$
\\
\bottomrule
\end{tabular}
\end{adjustbox}
\end{table*}

The analytical solution requires $15.0$,
$155.4$, and $98.4\,\mathrm{ms}$ for
$q_1$, $q_2$, and $q_3$, compared with $176.2$,
$457.9$, and $361.0\,\mathrm{ms}$ using
Crank--Nicolson.
The corresponding speedups are $11.8$,
$2.9$, and $3.7$ times.
The greatest speedup occurs for the lowest-order
system, for which the diagonalization cost
is reduced.
The analytical temporal evaluation provides
the temperatures at the observation times,
eliminating time stepping without
significantly changing the estimates.
By Corollary~\ref{cor:HN},
$\mathbf H_N(\hat{\mathbf a})$ is symmetric for every
$\hat{\mathbf a}$, so that the spectral decomposition
employed in the analytical solution is available at all
iterations.

\begin{figure}[!htbp]
\centering
\includegraphics[width=\columnwidth]{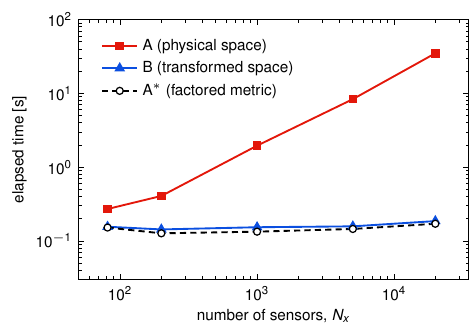}
\caption{Elapsed time in the iterative
Levenberg--Marquardt stage as a function of $N_x$,
for functionals A, B, and A$^*$.
Source $q_2$, step coefficient, and $N=32$.}
\label{fig:cost}
\end{figure}

Table~\ref{tab:cost} and Figure~\ref{fig:cost}
compare the computational costs of functionals A, B, and A$^*$.
Source $q_2$, the step coefficient,
and a fixed order $N=32$ are employed, with $N_t=101$,
varying only $N_x$.
The elapsed times correspond to the iterative stage and do not
include the initial processing of the measurements.
The residual vector dimension for A is
$N_xN_t$, whereas B and A$^*$ use
$NN_t$ components.

\begin{table*}[!tp]
\centering
\small
\caption{Computational cost of the iterative
Levenberg--Marquardt stage as a function of $N_x$:
filtered physical-domain functional A, transformed-domain
functional B, and factored evaluation A$^*$.
Source $q_2$, step coefficient, fixed $N=32$,
and $N_t=101$. The elapsed times do not include
the initial processing of the measurements.}
\label{tab:cost}
\setlength{\tabcolsep}{5pt}
\renewcommand{\arraystretch}{1.15}
\begin{adjustbox}{max width=\linewidth}
\begin{tabular}{rrrrrrrccc}
\toprule
&
\multicolumn{2}{c}{\textbf{Residual dimension}}
&
\multicolumn{3}{c}{\textbf{Iterative time [s]}}
&
\textbf{Ratio}
&
\textbf{It.}
&
\multicolumn{2}{c}{\textbf{NRMSE [\%]}}
\\
\cmidrule(lr){2-3}
\cmidrule(lr){4-6}
\cmidrule(lr){9-10}
$\boldsymbol{N_x}$
& \textbf{A} & \textbf{B}
& \textbf{A} & \textbf{B} & \textbf{A$^*$}
& \textbf{A/B} & \textbf{A/B}
& \textbf{A} & \textbf{B}
\\
\midrule
81
& 8181 & 3232
& 0.271 & 0.156 & 0.152
& 1.7 & 5/5 & 6.34 & 6.34
\\
200
& 20200 & 3232
& 0.407 & 0.143 & 0.127
& 2.8 & 5/5 & 5.81 & 5.81
\\
1000
& 101000 & 3232
& 1.959 & 0.153 & 0.134
& 12.8 & 5/5 & 5.68 & 5.68
\\
5000
& 505000 & 3232
& 8.366 & 0.158 & 0.145
& 53.0 & 5/5 & 5.71 & 5.71
\\
20000
& 2020000 & 3232
& 35.004 & 0.186 & 0.171
& 188.2 & 5/5 & 5.71 & 5.71
\\
\bottomrule
\end{tabular}
\end{adjustbox}
\end{table*}

The ratio between the elapsed times of A and B increases
from $1.7$, with $81$ sensors, to $12.8$,
with $1000$, and $188.2$, with $20000$.
In the latter case, the elapsed times are
$35.004$ and $0.186\,\mathrm{s}$.
The iterative cost of B remains approximately
constant because the order and the dimension
of the residual vector are kept fixed.
The initial processing, performed once before
the iterations, scales with
the number of measurements.

The A$^*$ evaluation exhibits elapsed times close
to those of B and reproduces A, with maximum
differences between the estimated coefficients,
divided by $\sqrt W$, on the order of
$10^{-15}\,\mathrm{s^{-1}}$.
Thus, the computational gain results from performing
the iterative operations in the reduced modal space
and can also be achieved while preserving the
metric of the filtered physical domain.
The agreement between the estimates from A and B
in the tests is a numerical result, whereas
the equivalence between A and A$^*$ follows from
Eq.~\eqref{eq:funcional_a_fatorado}.

In Table~\ref{tab:solver}, the elapsed times correspond
to the minimum value over five runs for
the analytical solution and three for Crank--Nicolson.
In the comparison of the functionals,
three runs per configuration are used for up to $5000$
sensors and one for $20000$.

\subsection{Comparison of truncation selection strategies}
\label{subsec:comparacao_truncamento}

To evaluate the accuracy and computational cost
of a priori selection relative to the
a posteriori criteria P1 and P2, the same
$270$ cases from Table~\ref{tab:estimates}
are considered.
Functional B, Levenberg--Marquardt,
and the analytical temporal solution are employed.
The reference configuration specifies $N_r=60$ for P2.
For this comparison, the $70$ modes mentioned in
Section~\ref{subsec:detalhes_inverso} are employed
and the range of candidate orders
is extended to $N_{\max}=70$, adopting
$N_r=70$ in the P2 evaluation.
This extension does not change the a priori selection,
whose maximum order remains $35$.

For each candidate order $n$, $n$
coefficients are estimated, and the final estimate
corresponds to the order $N_{\mathrm P}$ selected
by P1 or P2.
The NRMSE is calculated with $N=N_{\mathrm P}$.
The initialization from
Eq.~\eqref{eq:estimativa_inicial_coeficientes}
is applied by replacing $N$ with $n$.
Each order starts from the same initial estimate,
without reusing the coefficients from the previous
order. The search also allows the order
of minimum NRMSE to be determined, which is used only
as a reference because it requires knowledge
of $d(x)$.
Since $N_r=N_{\max}$, the P1 and P2 residuals
coincide at the last candidate order.

\begin{figure*}[!htbp]
\centering
\includegraphics[width=\textwidth]{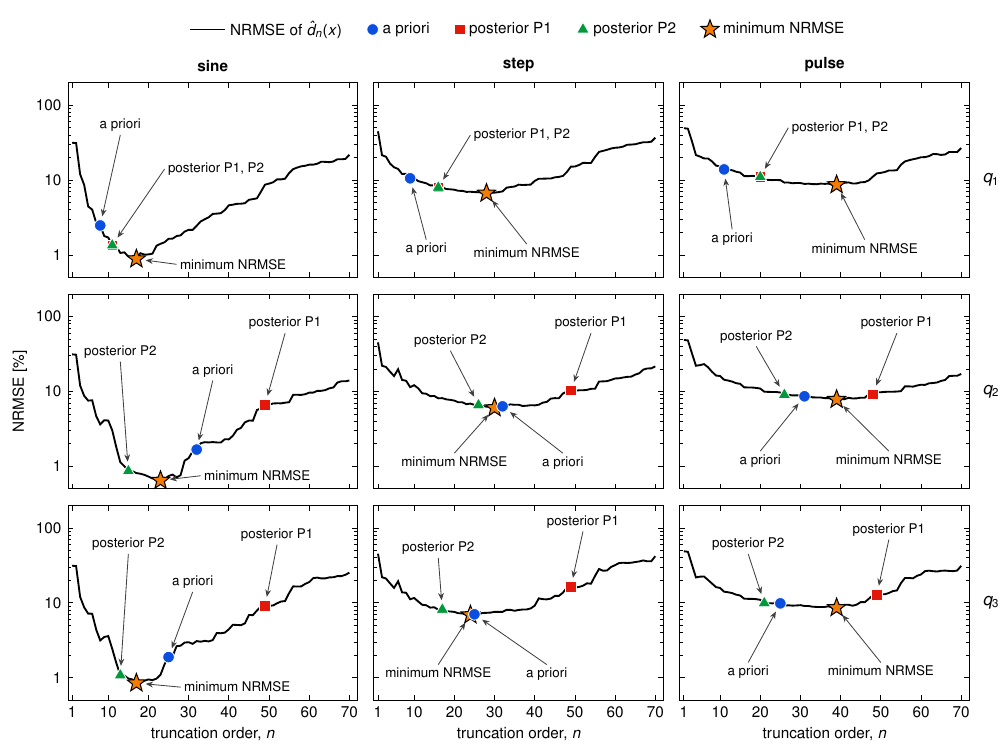}
\caption{NRMSE of $\hat d_n(x)$ as a function of
the order, for sources $q_1$, $q_2$, and $q_3$
in the rows and the sinusoidal, step, and pulse
profiles in the columns. The arrows indicate
the orders selected a priori, by P1 and P2,
and by the minimum NRMSE.
Results for $\sigma=0.10\,^\circ\mathrm C$
and random seed equal to one.}
\label{fig:error_order}
\end{figure*}

\begin{figure*}[!htbp]
\centering
\includegraphics[width=\textwidth]{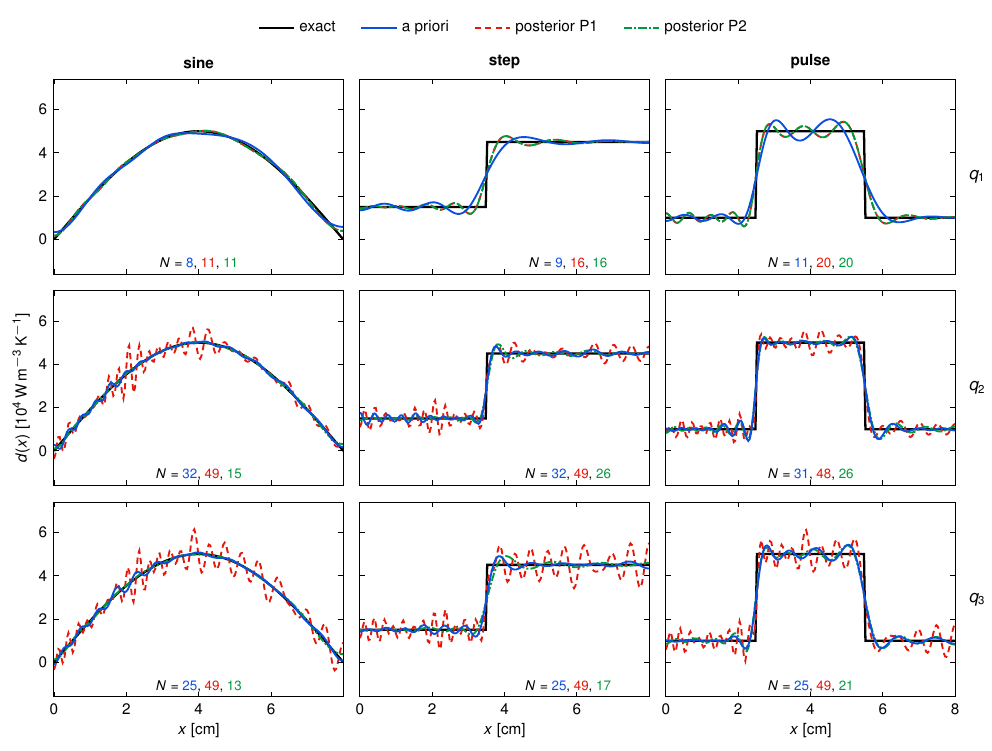}
\caption{Reconstructions $\hat d_N(x)$ using
a priori and a posteriori selection, for
sources $q_1$, $q_2$, and $q_3$ in the rows
and the sinusoidal, step, and pulse profiles
in the columns.
For the a posteriori criteria, $N=N_{\mathrm P}$.
The orders are indicated by the curve colors.
Results for $\sigma=0.10\,^\circ\mathrm C$
and random seed equal to one.}
\label{fig:posterior_estimates}
\end{figure*}

Figures~\ref{fig:error_order}
and~\ref{fig:posterior_estimates} use
$81\times101$ measurements, $\sigma=0.10\,^\circ\mathrm C$,
and the realization with random seed equal to one.

Figure~\ref{fig:error_order} presents
the error of $\hat d_n(x)$ as a function of the order,
together with the selections of each strategy.
For the step profile with $q_2$ and $q_3$, the
a priori selection yields $N=32$ and $25$, close
to the orders of minimum error, $30$ and $24$.
P1 selects $49$ terms in both cases,
with errors of $10.25\%$ and $16.33\%$,
compared with $6.34\%$ and $7.04\%$ for the
a priori selection.
P2 selects $26$ and $17$ terms, with errors
of $6.59\%$ and $8.15\%$.

The a priori selection does not always
coincide with the order of minimum error.
With $q_1$, it retains fewer terms than those
corresponding to the minimum; for the sinusoidal
profile with $q_2$ and $q_3$, it retains more terms.
In the realization presented, the a posteriori
criteria improve the reconstructions with $q_1$,
and P2 yields the lowest error
among the three strategies for the sinusoidal
profile with the zonal sources.

Figure~\ref{fig:posterior_estimates} shows
the reconstructed functions.
The high orders selected by P1 with $q_2$ and $q_3$
introduce short-wavelength oscillations,
whereas the lower orders selected by P2
produce greater smoothing of the transitions.
With $q_1$, the a posteriori estimates
represent the jumps with steeper transitions.

Table~\ref{tab:comparison_summary} summarizes
the results of the $270$ cases.
The a priori selection yields a mean NRMSE
of $6.99\%$, lower than the $14.38\%$
of P1 and $12.34\%$ of P2.
The mean error using the order of minimum NRMSE
is $5.35\%$, whereas using
$70$ terms in all cases results
in $31.12\%$.
Therefore, extending the representation
to the available limit can significantly
deteriorate the estimates.

\begin{table*}[!tp]
\centering
\small
\caption{Comparison of the strategies over $270$
cases, with $N_{\max}=N_r=70$:
mean selected order, mean NRMSE,
number of cases in which the discrepancy condition
was not satisfied, mean number of estimations,
and total computational time relative to a priori selection.
For P1 and P2, the order corresponds to
$N_{\mathrm P}$.
The minimum NRMSE reference uses the exact function.}
\label{tab:comparison_summary}
\setlength{\tabcolsep}{5pt}
\renewcommand{\arraystretch}{1.15}
\begin{adjustbox}{max width=\linewidth}
\begin{tabular}{lcccccccc}
\toprule
& &
\multicolumn{4}{c}{\textbf{NRMSE [\%]}}
& \textbf{Not}
&
& \textbf{Relative}
\\
\cmidrule(lr){3-6}
\textbf{Criterion}
& \textbf{Mean order}
& \textbf{Sinusoidal}
& \textbf{Step}
& \textbf{Pulse}
& \textbf{All}
& \textbf{satisfied}
& \textbf{Estimations}
& \textbf{time}
\\
\midrule
A priori
& 21.3 & 1.82 & 8.09 & 11.06 & 6.99
& 0 & 1 & 1
\\
Posterior P1
& 47.5 & 9.94 & 17.95 & 15.23 & 14.38
& 90 & 49.2 & 182
\\
Posterior P2
& 35.9 & 7.54 & 15.25 & 14.25 & 12.34
& 90 & 38.6 & 162
\\
\midrule
Minimum NRMSE
& 28.1 & 0.88 & 6.60 & 8.58 & 5.35
& -- & -- & --
\\
\bottomrule
\end{tabular}
\end{adjustbox}
\end{table*}

The a priori selection satisfies the discrepancy
condition in all cases.
P1 and P2 reach the prescribed level in $180$
cases each, with agreement regarding satisfaction
of the condition in $268$ of the
$270$ cases.
In the $90$ cases without a crossing for each
criterion,
Eq.~\eqref{eq:selecao_posterior_minimo}
selects orders between $59$ and $70$ for P1
and between $50$ and $70$ for P2, with means
of $64.8$ and $62.0$, respectively.
These orders lie beyond the region
of minimum error in many cases.

The a posteriori residuals approach
the noise level with a small slope,
making the crossing more sensitive to
perturbations.
For the number of measurements considered,
the mean of the squared Gaussian errors
has a relative standard deviation of
approximately $1.6\%$ (Lemma~\ref{lem:disc}).
This fluctuation, together with model
and estimation errors, contributes to
the variation in the a posteriori selected orders.
When the discrepancy condition is satisfied
only at high orders, the reconstruction of
$d(x)$ is less accurate.

Propositions~\ref{prop:apriori}
and~\ref{prop:aposteriori} identify the origin of this
behavior. For $81\times101$ measurements, each retained mode
removes $1.23\%$ of $\sigma^2$ from the a priori residual,
but only $0.012\%$ from the a posteriori residuals, whereas
the fluctuation of the discrepancy is $1.56\%$ of
$\sigma^2$. By Corollary~\ref{cor:crossing} with $\tau=1$,
even an exact model satisfies the discrepancy condition with
probabilities of only $0.53$, $0.61$, $0.65$, and $0.71$ at
$n=10$, $35$, $50$, and $70$. The occurrence of a crossing
is thus decided essentially by whether the sample variance
of the error realization lies below $\sigma^2$.

In the simulations, the standardized error field
$\varepsilon_{sj}/\sigma$ depends only on the realization
index, and not on the source, the profile, or $\sigma$, so
that the $270$ cases contain ten independent error fields.
Their sample variance is below $\sigma^2$ in six of them,
and the occurrence of a crossing coincides with this event
in $244$ (P1) and $242$ (P2) of the $270$ cases; all
exceptions correspond to realizations whose sample variance
is within $1\%$ of $\sigma^2$. This accounts for the $180$
cases in which each criterion reaches the prescribed level
and for their agreement in $268$ cases. The cases without a
crossing therefore reflect the absence of a tolerance
factor, and not the quality of the estimates.

The a posteriori criteria were therefore also evaluated with
Eq.~\eqref{eq:selecao_posterior_tau}. Since the search over
$n=1,\ldots,70$ provides the residuals and the errors at
every order, the rule
$\sigma_{\mathrm P,n}^2\le\tau\sigma^2$ can be applied to
the same estimates. Table~\ref{tab:tolerance} summarizes the
results for $\tau=1+3\sqrt{2/(N_xN_t)}\approx1.05$
(Remark~\ref{rem:tolerance}). No case remains without a
crossing, for P1 and for P2, which already holds for
$\tau\ge1.03$. The mean NRMSE becomes $7.45\%$ for P1 and
$7.74\%$ for P2, compared with $14.38\%$ and $12.34\%$ for
$\tau=1$ and $6.99\%$ for the a priori selection, and the
mean number of estimations decreases from $49.2$ and $38.6$
to $30.6$ and $13.5$, with computational times $24$ and $5$
times that of the a priori selection.

By profile, P1 yields $3.37\%$ (sinusoidal), $8.79\%$
(step), and $10.20\%$ (pulse), compared with $1.82\%$,
$8.09\%$, and $11.06\%$ for the a priori selection. With
$q_1$, P1 is more accurate than the a priori selection for
the three profiles, with means over the three noise levels
of $1.96\%$, $9.90\%$, and $12.79\%$, against $2.39\%$,
$10.75\%$, and $14.10\%$. This result agrees with
Remark~\ref{rem:scaling}: for a spatially uniform source,
the order selected a priori follows the smooth temperature
field and remains below the order of minimum error. With the
zonal sources, P1 selects orders between $22$ and $56$ and
is less accurate for the sinusoidal and step profiles. For
$\tau$ between $1.02$ and $1.20$, the mean error remains
between $7.1\%$ and $8.3\%$ for P1 and between $7.0\%$ and
$9.5\%$ for P2, never below the value of the a priori
selection.

Most of the difference in accuracy reported in
Table~\ref{tab:comparison_summary} is therefore due to the
absence of a tolerance factor in
Eq.~\eqref{eq:selecao_posterior}. With $\tau>1$, the
a priori selection retains a lower mean error, by $0.5$ to
$0.8$ percentage point, and requires a single estimation.

\begin{table*}[!tp]
\centering
\small
\caption{Effect of the tolerance factor in
Eq.~\eqref{eq:selecao_posterior_tau} on the a posteriori
criteria over the $270$ cases, with
$\tau=1+3\sqrt{2/(N_xN_t)}\approx1.05$: mean selected order,
mean NRMSE, number of cases in which the condition was not
satisfied, mean number of estimations, and total
computational time relative to a priori selection.}
\label{tab:tolerance}
\setlength{\tabcolsep}{5pt}
\renewcommand{\arraystretch}{1.15}
\begin{adjustbox}{max width=\linewidth}
\begin{tabular}{lccccccccc}
\toprule
& & &
\multicolumn{4}{c}{\textbf{NRMSE [\%]}}
& \textbf{Not}
&
& \textbf{Relative}
\\
\cmidrule(lr){4-7}
\textbf{Criterion}
& $\boldsymbol\tau$
& \textbf{Mean order}
& \textbf{Sinusoidal}
& \textbf{Step}
& \textbf{Pulse}
& \textbf{All}
& \textbf{satisfied}
& \textbf{Estimations}
& \textbf{time}
\\
\midrule
A priori & -- & 21.3 & 1.82 & 8.09 & 11.06 & 6.99
& 0 & 1 & 1
\\
\midrule
Posterior P1 & 1 & 47.5 & 9.94 & 17.95 & 15.23 & 14.38
& 90 & 49.2 & 182
\\
& 1.05 & 30.6 & 3.37 & 8.79 & 10.20 & 7.45
& 0 & 30.6 & 24
\\
\midrule
Posterior P2 & 1 & 35.9 & 7.54 & 15.25 & 14.25 & 12.34
& 90 & 38.6 & 162
\\
& 1.05 & 13.5 & 2.12 & 9.44 & 11.64 & 7.74
& 0 & 13.5 & 5
\\
\bottomrule
\end{tabular}
\end{adjustbox}
\end{table*}

Table~\ref{tab:comparison_cases} details
the behavior for
$\sigma=0.10\,^\circ\mathrm C$.
For the step profile with $q_1$, the a priori
selection yields $N=9$ in all ten realizations,
whereas P1 ranges from $11$ to $63$.
For the step profile with $q_3$, the standard
deviations of the NRMSE are $12.23\%$ and $14.27\%$
for P1 and P2, compared with $0.17\%$ for
the a priori selection.
In these nine combinations of source and profile,
both a posteriori criteria exhibit higher
mean errors than the a priori selection,
although they may provide better estimates
in individual realizations.

\begin{table*}[!tp]
\centering
\small
\caption{Comparison for
$\sigma=0.10\,^\circ\mathrm C$,
with $N_{\max}=N_r=70$:
ranges of selected orders, mean NRMSE,
and population standard deviation over ten realizations.
Rel.~time denotes the total computational time
relative to a priori selection.
For the minimum NRMSE reference, the range
of orders and mean error are presented.}
\label{tab:comparison_cases}
\setlength{\tabcolsep}{5pt}
\renewcommand{\arraystretch}{1.15}
\begin{adjustbox}{max width=\linewidth}
\begin{tabular}{llcccccccccc}
\toprule
& &
\multicolumn{2}{c}{\textbf{A priori}}
&
\multicolumn{3}{c}{\textbf{Posterior P1}}
&
\multicolumn{3}{c}{\textbf{Posterior P2}}
&
\multicolumn{2}{c}{\textbf{Minimum NRMSE}}
\\
\cmidrule(lr){3-4}
\cmidrule(lr){5-7}
\cmidrule(lr){8-10}
\cmidrule(lr){11-12}
\textbf{Source}
& $\boldsymbol{d(x)}$
& $N$
& \textbf{NRMSE [\%]}
& $N_{\mathrm{P1}}$
& \textbf{NRMSE [\%]}
& \textbf{Rel.~time}
& $N_{\mathrm{P2}}$
& \textbf{NRMSE [\%]}
& \textbf{Rel.~time}
& \textbf{Order}
& \textbf{NRMSE [\%]}
\\
\midrule
$q_1$ & sinusoidal
& 7--8 & 2.56 $\pm$ 0.10
& 9--65 & 6.77 $\pm$ 7.97 & 328
& 9--66 & 6.86 $\pm$ 8.05 & 351
& 14--20 & 0.89
\\
& step
& 9 & 10.56 $\pm$ 0.01
& 11--63 & 14.43 $\pm$ 9.64 & 255
& 12--63 & 14.26 $\pm$ 9.47 & 258
& 20--33 & 6.85
\\
& pulse
& 11 & 13.88 $\pm$ 0.01
& 16--66 & 14.09 $\pm$ 5.89 & 266
& 16--66 & 14.10 $\pm$ 5.89 & 270
& 26--39 & 8.82
\\
\midrule
$q_2$ & sinusoidal
& 31--32 & 1.41 $\pm$ 0.37
& 43--65 & 6.83 $\pm$ 1.43 & 98
& 13--65 & 3.68 $\pm$ 3.63 & 60
& 19--23 & 0.72
\\
& step
& 31--32 & 6.19 $\pm$ 0.18
& 43--65 & 11.55 $\pm$ 2.23 & 101
& 18--65 & 9.35 $\pm$ 3.65 & 69
& 28--42 & 5.92
\\
& pulse
& 31 & 8.55 $\pm$ 0.09
& 42--64 & 9.82 $\pm$ 0.93 & 92
& 21--65 & 9.69 $\pm$ 1.37 & 64
& 39--46 & 7.74
\\
\midrule
$q_3$ & sinusoidal
& 24--25 & 1.82 $\pm$ 0.24
& 36--70 & 14.45 $\pm$ 6.69 & 150
& 12--69 & 9.10 $\pm$ 10.73 & 96
& 14--20 & 0.94
\\
& step
& 24--26 & 7.06 $\pm$ 0.17
& 36--70 & 24.73 $\pm$ 12.23 & 141
& 12--69 & 18.29 $\pm$ 14.27 & 90
& 23--30 & 6.93
\\
& pulse
& 24--25 & 9.92 $\pm$ 0.08
& 36--70 & 18.84 $\pm$ 8.65 & 147
& 16--69 & 15.81 $\pm$ 8.49 & 97
& 30--39 & 8.86
\\
\bottomrule
\end{tabular}
\end{adjustbox}
\end{table*}

For the datasets in which the discrepancy
condition is satisfied, P1 and P2 exhibit mean
errors of $10.24\%$ and $7.35\%$, respectively,
whereas the a priori selection yields
approximately $7.07\%$ in the corresponding cases.
Without a crossing, the mean errors increase
to $22.65\%$ and $22.33\%$,
compared with approximately $6.82\%$ for
the a priori selection.
Thus, P2 approaches the mean accuracy
of the a priori selection when it satisfies
the discrepancy condition, but exhibits
substantial deterioration when the selection depends on
Eq.~\eqref{eq:selecao_posterior_minimo}.

Regarding computational cost, P1 and P2 require,
on average, $49.2$ and $38.6$ estimations,
with elapsed times $182$ and $162$ times greater
than that of the a priori selection.
The elapsed times include the discrepancy
calculation and one estimation in the a priori
selection, and successive estimations and residual
evaluations in the a posteriori criteria.
In cases without a crossing, $70$
estimations are performed.
The time ratios exceed the mean number
of estimations because the cost of
each fit increases with the order.

The results correspond to the discrepancy
level without a tolerance factor and
independent estimations for each order.
Reusing the coefficients from the previous
order may reduce the iterative cost of the
successive inverse problems.
Under the conditions examined, the a priori
selection exhibits lower mean error,
less variation in the selected order across realizations,
and lower computational cost.
Two observations complement this comparison.
First, with a tolerance factor, the a posteriori criteria
approach the accuracy of the a priori selection, which
retains an advantage of $0.5$--$0.8$ percentage point in
mean error (Table~\ref{tab:tolerance}) and a single
estimation per case. Second, since the order
selected a priori follows the spatial content of the
temperature field (Remark~\ref{rem:scaling}), spatially
uniform excitation leads to orders below that of minimum
error: for the step with $q_1$ and
$\sigma=0.10\,^\circ\mathrm C$, $N=9$ is selected, whereas
the minimum NRMSE occurs between orders $20$ and $33$, with
mean errors of $10.56\%$ and $6.85\%$, respectively
(Table~\ref{tab:comparison_cases}). In this configuration,
criterion P1 with $\tau>1$ provides the most accurate
reconstructions.

\section{Conclusions}\label{sec:conclusions}

This work investigated the identification of spatially varying
thermal dissipation coefficients in heterogeneous media
using integral transforms and nonlinear least squares.
Representing the unknown coefficient in the same spectral
basis employed for the temperature allowed modal coupling
to be expressed through precomputed matrices,
while sensitivities obtained from the transformed system
were used in the Gauss--Newton and Levenberg--Marquardt
procedures. This formulation made it possible to examine
the influence of spectral truncation on both the
representation of the solution and the number of
estimated parameters.

\textit{A priori} and \textit{a posteriori} strategies
for selecting the truncation order based on the
discrepancy principle were investigated. The
\textit{a priori} selection, performed directly
on the measured temperatures, allowed the order
to be determined before optimization, avoiding
the solution of successive inverse problems.
Across the $270$ cases evaluated, this strategy
achieved a mean NRMSE of $6.99\%$ with one estimation
per case.
With a tolerance factor of approximately $1.05$ on the
discrepancy level, the a posteriori criteria P1 and P2
satisfied the condition in all cases, with mean errors of
$7.45\%$ and $7.74\%$ and computational times $24$ and $5$
times greater than that of the a priori selection. Without
the tolerance factor, the condition was not reached in $90$
cases for each criterion, and the mean errors and relative
times increased to $14.38\%$ and $12.34\%$ and to $182$ and
$162$, respectively. The analysis of the residuals explains
this behavior: each retained mode removes a fraction
$1/N_x$ of the noise variance from the a priori residual,
against $1/(N_xN_t)$ for the a posteriori residuals, which
is below the statistical fluctuation of the discrepancy.
For spatially uniform excitation, the order selected
\textit{a priori} followed the smooth temperature field and
remained below the order of minimum error, and criterion P1
with the tolerance factor provided the most accurate
reconstructions.

The reconstructions showed that accuracy depends
on the spatial regularity of the unknown coefficient,
the thermal excitation, and the measurement resolution.
Smooth profiles exhibited lower errors, whereas
discontinuities were accompanied by smoothing
and oscillations associated with spectral truncation.
The truncation error of the expansion decays as
$N^{-1/2}$ for discontinuous coefficients and as $N^{-3/2}$
for the smooth profile considered, and it accounted for most
of the reconstruction error of the step and pulse profiles
at the selected orders.
Spatially localized sources favored the identification
of higher-frequency components, particularly with
adequate spatial resolution. The results also showed
that increasing the truncation order does not
ensure greater accuracy and may amplify
the effects of noise on the reconstructed function.

The comparison of the least-squares functionals
demonstrated the computational advantages of fitting
in the transformed domain. For $20000$ spatial sensors,
the time required for the iterative stage was reduced
from $35.004$~s in the filtered physical domain
to $0.186$~s in the transformed domain.
Furthermore, the factored modal evaluation
reproduced the filtered physical-domain functional,
preserving its metric at a computational cost
close to that of the transformed formulation.
The analytical temporal solution also provided
computational speedups ranging from $2.9$ to
$11.8$ times relative to Crank--Nicolson integration,
without significant differences in the
resulting estimates.

Overall, the results indicate that prior truncation
selection, combined with spectral parameterization
and residual evaluation in modal space, provides
an efficient alternative for the nonlinear
identification of spatially varying thermal
dissipation coefficients. The methodology allows
the dimensionality of the inverse problem to be
controlled before estimation, reduces computational
cost, and yields reconstructions with satisfactory
accuracy in the presence of noise.
Application to experimental data and extension
to multidimensional problems represent
possible directions for future work.

\ifproofsinline\else
\appendix
\section{\corr{Proofs}}
\label{app:proofs}
\AppendixProof{Lemma~\ref{lem:eig}}{\ProofLemEig}
\AppendixProof{Proposition~\ref{prop:H}}{\ProofPropH}
\AppendixProof{Proposition~\ref{prop:decomp}}{\ProofPropDecomp}
\AppendixProof{Theorem~\ref{thm:step}}{\ProofThmStep}
\AppendixProof{Proposition~\ref{prop:neumann}}{\ProofPropNeumann}
\AppendixProof{Corollary~\ref{cor:HN}}{\ProofCorHN}
\AppendixProof{Lemma~\ref{lem:disc}}{\ProofLemDisc}
\AppendixProof{Proposition~\ref{prop:apriori}}{\ProofPropApriori}
\AppendixProof{Proposition~\ref{prop:aposteriori}}{\ProofPropAposteriori}
\AppendixProof{Corollary~\ref{cor:crossing}}{\ProofCorCrossing}
\AppendixProof{Theorem~\ref{thm:ident}}{\ProofThmIdent}
\fi

\selectlanguage{english}
\section*{Acknowledgements}

This study was financed in part by CAPES - Coordenação de Aperfeiçoamento de Pessoal de Nível Superior - Brazil, Finance Code 001. The authors also thank the other sponsoring agencies, CNPq - Conselho Nacional de Desenvolvimento Científico e Tecnológico - Brazil, FAPERJ - Fundação Carlos Chagas Filho de Amparo à Pesquisa do Estado do Rio de Janeiro, FAPESP - Fundação de Amparo à Pesquisa do Estado de São Paulo, SENAI CIMATEC, and IFNMG - Instituto Federal do Norte de Minas Gerais - Brazil.

\justifying
\bibliographystyle{ieeetr}
\bibliography{references}

\end{document}